\documentclass[11pt,fleqn]{article}
\usepackage{graphicx}
\usepackage{amsmath}
\usepackage{amsthm}
\usepackage{amssymb}
\usepackage{amsfonts}
\usepackage{epsfig}
\usepackage{hyperref}

\usepackage{epstopdf}

\usepackage{color}

\def\Im {\mathop{\rm Im}\nolimits}
\def\arg {\mathop{\rm arg}\nolimits}
\def\Re {\mathop{\rm Re}\nolimits}
\def\Ai {{\rm Ai}}
\def\Bi {{\rm Bi}}

\newtheorem{rem}{Remark}
\newtheorem{lem}{Lemma}
\newtheorem{thm}{Theorem}
\newtheorem{cor}{Corollary}

\numberwithin{equation}{section}

\title{Asymptotics of the Charlier polynomials  via difference equation  methods}

\author{Xiao-Min Huang$^{\,1,\,2 }$, Yu Lin$^{\,3}$  and  Yu-Qiu Zhao$^{\,4}$}

\date{
\hbox{\small \emph{$^1$  School of Applied Mathematics, Guangdong University of Technology, Guangzhou,  China}}
\hbox{\small \emph{$^2$ Department of Mathematics, City University of Hong Kong, Hong Kong }}
 \hbox{\small
\emph{$^3$ Department of Mathematics, South China University of Technology, Guangzhou, China}
 }
 \hbox{\small
\emph{$^4$ Department of Mathematics, Sun Yat-sen University, Guangzhou,
China}}}

\begin{document}

\maketitle
\begin{abstract}
We derive  uniform and non-uniform  asymptotics of the Charlier polynomials by using difference equation methods alone. The   Charlier polynomials are special in that they
   do not fit into the framework of the turning point theory,  despite the fact that they are crucial  in the Askey scheme. In this paper, asymptotic approximations are obtained respectively in the outside
region, an intermediate region, and near the turning points.
In particular, we  obtain uniform asymptotic approximation at a pair of coalescing turning points with the aid of a local transformation.
  We also give a
 uniform approximation at the origin by  applying the method of dominant balance and several matching techniques.
\end{abstract}

\vspace{.4cm}



\section{Introduction and statement of results}

The Charlier polynomials $C^{(a)}_n(x)$  are   discrete orthogonal polynomials such that
\begin{equation}\label{Charlier-orthogonal}
\sum^\infty_{k=0} C^{(a)}_n(k)   C^{(a)}_m(k) \frac {a^k}{k!}=e^a a^n n! \delta_{m n},~~m, ~n=0,1,\cdots,
\end{equation}with parameter  $a>0$.
An explicit expression for the polynomials is
\begin{equation}\label{Charlier-explicit}
C^{(a)}_n(x)=\sum^n_{k=0} \begin{pmatrix}  n \\
                                                                            k \\
                                                                          \end{pmatrix}
\begin{pmatrix}
                                                                            x \\
                                                                            k \\
                                                                          \end{pmatrix}
k!(-a)^{n-k};
\end{equation}cf. Szeg\H{o} \cite[pp.34-35]{szego-book}.  Here the notation $C^{(a)}_n(x)$ refers to the monic polynomial,  as in Bo and Wong \cite{bo-wong1994}.  Other notations are also used to stand for the Charlier polynomials. For example, $C_n(x, a)=(-a)^{-n} C_n^{(a)} (x)$; cf. \cite[Ch.18]{nist}. The family of Charlier polynomials occupies a crucial position in the Askey scheme, as illustrated in \cite[Fig.18.21.1]{nist}.

In 1994,  Bo and Wong \cite{bo-wong1994} considered the uniform asymptotic expansions for $C^{(a)}_n(n\beta)$ as $n\to\infty$. The uniformity is for $\beta$ in compact subsets of the real interval $(0, +\infty)$. Integral methods are   used to obtain the asymptotic formulas, starting with the
 generating function
\begin{equation}\label{charlier-generating}
e^{-aw} (1+w)^x=\sum^\infty_{n=0} C^{(a)}_n(x) \frac {w^n}{n!}.
\end{equation}
The uniform interval of   \cite{bo-wong1994}   covers the seven regions in a 1998 paper \cite{goh1998} of Goh. In these regions,
Plancherel-Rotach asymptotics are obtained for the Charlier
polynomials,  using the integral methods as well.

It is readily seen that   \eqref{Charlier-explicit}  follows from  \eqref{charlier-generating}; cf.\;\cite{szego-book}. Also, Darboux's method and the steepest descent method can be applied  to derive non-uniform asymptotic approximations.  Earlier in 1985, an asymptotic formula for
$C_n^{(a)}(x)$ when $x < 0$ has been obtained   by Maejima and Van Assche \cite{maejima-vanassche1985}, using
probabilistic arguments.

In \cite{bo-wong1994},
Bo and Wong commented that {\it in regard to the asymptotics of the Charlier polynomials, not much is known in the literature.} Since then, a quarter century sees  new observations made and quite a number of novel tools developed.
For instance, in 2001, Dunster  \cite{dunster2001} made use of the connection  with Laguerre  polynomials, namely,
$C_n^{(a)}(x)=n! L_n^{(x-n)}(a)$,  considered a  differential equation with respect to the parameter $a$, and obtained asymptotic formulas, uniformly in $x$ or $a$.

The Riemann-Hilbert approach is a powerful new tool in asymptotic analysis; see Deift \cite{deift1999}.  For discrete orthogonal polynomials,   pioneering work has been done
 by Baik et al.\;\cite{baik-et-al2007}, published  in 2007.  The method of Baik et al. is applied by  Ou and Wong \cite{ou-wong2010}, with modifications, to obtain uniform asymptotic approximations for $C_n^{(a)} (z)$ with non-rescaled complex variable $z$. A significant feature of \cite{ou-wong2010} is the global uniformity. Asymptotic expansions are derived in three  regions that cover the whole complex plane.
 Yet formulas with  explicit   leading coefficients are not written down.
   It is worth mentioning that
 the Charlier polynomials demonstrate a uniform  equilibrium measure,
that is, the asymptotic zero distribution of  $C_n^{(a)}(ny)$ has a constant density $1$, supported on $y\in [0,1]$; cf. \cite{kuijlaars-vanassche1999}.

  Thus, up to now, quite a lot  of facts have been known about the Charlier polynomials, including  generating functions, differential equations, recurrence relations, and the orthogonal measure.  Uniform and non-uniform asymptotics are derived. Still, the global uniform asymptotic approximations need to be clarified, and the turning point asymptotics of the polynomials are of great interest.  The polynomials can, in a sense,  serve as a touchstone    for new tools and techniques developed.

  The main focus of the present investigation will be on difference equation methods.

  In a  review \cite{askey1985} of Chihara's book \cite{chihara-book}, Richard Askey commented on the reason for the renewal of interest in orthogonal polynomials  in the 70-80s, that {\it General orthogonal polynomials are primarily interesting because of their 3-term recurrence relation}.

The three-term  recurrence formula
for  Charlier polynomials is
\begin{equation}\label{charlier-difference-eq}
xC^{(a)}_n(x)=  C^{(a)}_{n+1}(x)+(n+a)  C^{(a)}_n(x)+an   C^{(a)}_{n-1}(x),~~n=0,1,\cdots,
\end{equation}with fixed $a>0$, and initial data $C^{(a)}_{-1}=0$ and $C^{(a)}_{0}=1$; cf. \cite{bo-wong1994, chihara-book}.

\subsection{Non-oscillatory regions and a neighborhood of the origin}

A non-oscillatory region   considered here is described as
\begin{equation*}
x=ny, ~~~~~|y-[0,1]|>r
\end{equation*}
  for  an arbitrary  positive constant $r$. Here $|y-[0,1]|$ stands for the distance between the point $y$ and the line segment $[0, 1]$.
The natural re-scaling   $x=ny$ aims to normalize the limiting oscillatory interval to $y\in [0, 1]$; see Kuijlaars and Van Assche \cite[Sec.4.5]{kuijlaars-vanassche1999}.
There are several ways to derive  the asymptotics in these unbounded   regions via difference equation.
In what follows we take
 arguments from Wang and Wong \cite{wang-wong2012}, and   in sprit  similar to  Van Assche and Geronimo \cite{vanassche-geronimo1989}.

Denote the coefficients in \eqref{charlier-difference-eq} as  $a_n=n+a$ and $b_n= an$, and introduce
\begin{equation}\label{approx-out-Charlier-monic}
C^{(a)}_n(x)= \prod ^n_{k=1} w_k(x).
\end{equation}We see that $w_1(x)=x-a$, and
\begin{equation}\label{w-k-Charlier-first-order}
w_{k}(x)=x-a_{k-1}-\frac {b_{k-1}} {w_{k-1}(x)},~~k=2,3,\cdots.
\end{equation}
It is   verified in  Lemma \ref{lemma-w-k-estimate},  that
\begin{equation}\label{w-k-Charlier}
w_k=(x-k)\left\{ 1+\frac {(1-a) x-k}{(x-k)^2} +O\left (\frac 1 {n^2}\right )\right\},
\end{equation}in which the error term is uniform in $k=1,2,\cdots, n$, and in $|y- [0, 1]|>r$ for arbitrary positive $r$; see Section \ref{subsec:non-oscillatory-asymptotics} for full details of the proof.

Now substituting \eqref{w-k-Charlier} into \eqref{approx-out-Charlier-monic}, and  making use of the trapezoidal rule, we obtain the asymptotic formula
\begin{equation}\label{approx-out-Charlier-formula}
  C^{(a)}_n(n y)  =    n^{n} \sqrt{\frac y {y-1}}  \exp\left (-\frac {a}{y-1}\right ) \exp\left\{n  \left [   y  \log\frac {  y }{  y-1} -1\right ] \right\} (y-1)^n\left [1+O\left(\frac 1 n\right )\right ],
\end{equation}holding uniformly for large $n$ and  $y$ keeping a constant distance from  $[0, 1]$. The logarithm  and square root   take principal branches.


Throughout this paper, we stick to our theme, using difference equation methods alone.
For instance,
 Section \ref{subsec:origin}  will be devoted to a uniform asymptotic analysis of $C^{(a)}_n(x)$ in a neighborhood of $x=0$. To this end, we apply the method of  dominant balance to the difference equation \eqref{charlier-difference-eq}, to obtain a subdominant
asymptotic solution $C_s(n,x)$, and a dominant asymptotic solution $C_d(n,x)$.
The latter is then extended   uniformly to a domain $|y|\leq r <1$, with $x=ny$. The uniform asymptotic approximation for $C^{(a)}_n(x)$ in the domain of uniformity is then determined by matching the outer asymptotics \eqref{approx-out-Charlier-formula},  as stated in  Lemma \ref{prop-approx-at-0}.

A combination of   formulas \eqref{approx-out-Charlier-formula} with \eqref{Charlier-at-0}
  yields
the following uniform asymptotic approximation:

\begin{thm}\label{thm-approx-not-at-1} For arbitrary positive constant $\delta$,  it   holds
\begin{equation}\label{approx-not-at-1}
C^{(a)}_n(ny)  = (-1)^n e^{ \frac {a }{1-y }}\frac {\Gamma\left (n-ny\right )}{\Gamma(-ny)}\left[1+O\left(\frac 1 n\right)\right]+ \varepsilon_1,~~~\mbox{for}~~|y-1|> \delta,
\end{equation}
where
$\varepsilon_1$   is exponentially small
 as $\left |y-[0,1]\right |> r$ for arbitrary   positive $r$.
 \end{thm}

\begin{figure}[t]
\centering
\includegraphics[width=7.8   cm]{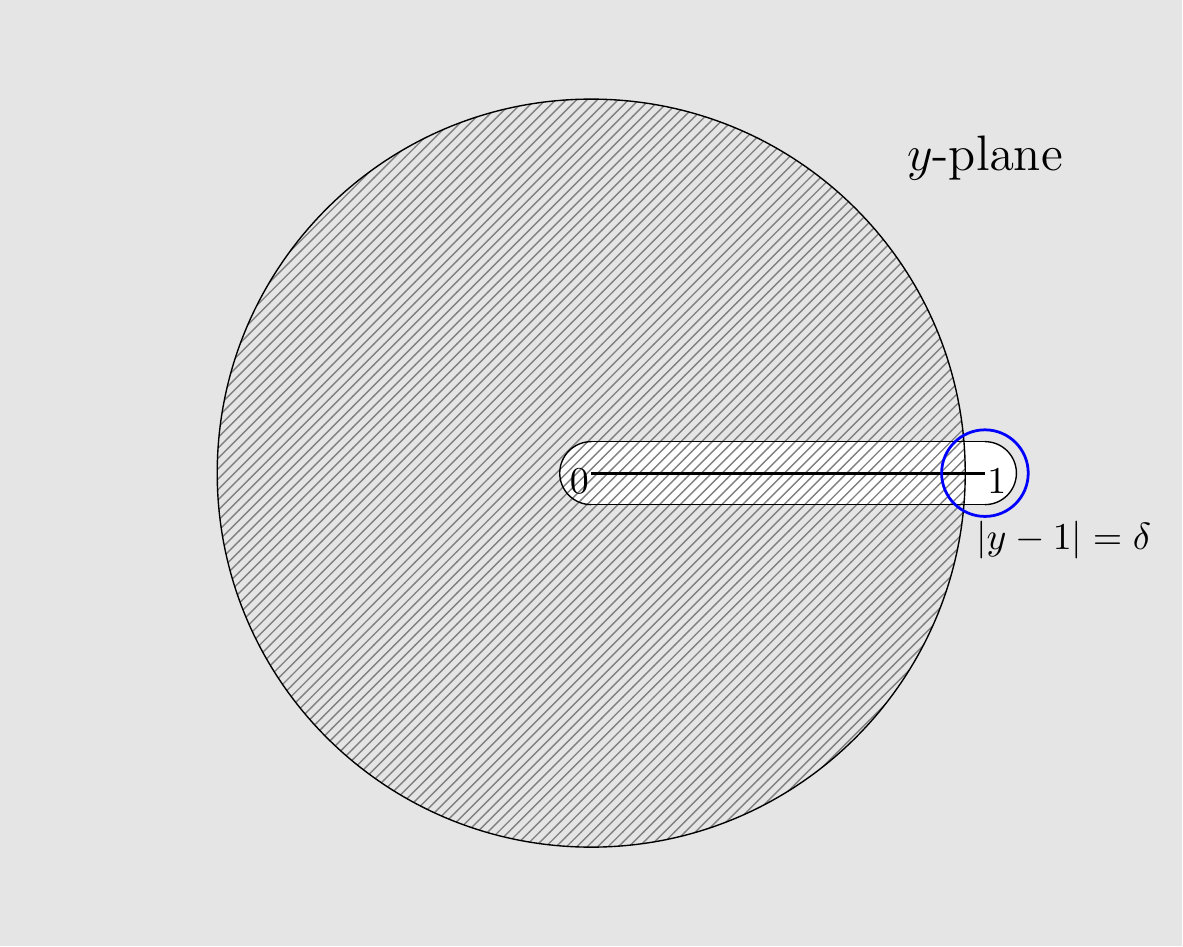}\hskip 1.1cm \includegraphics[width=6.8    cm]{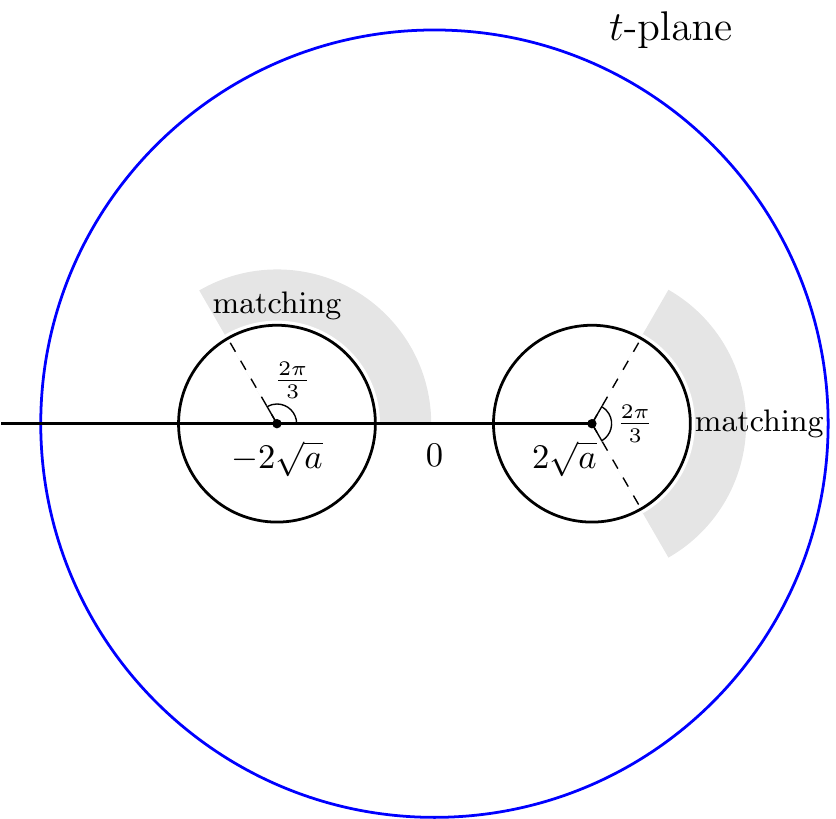}
\caption{\footnotesize Left: Domains of validity. $D_1$ is the shaded disc $|y|<1-\delta'$, $D_2$ denotes the outer region  $|y-[0, 1]| >\delta''$ in gray, then $D_1\cup D_2$ covers  $|y-1|\geq \delta$, with small positive $\delta$ and appropriately chosen $\delta'$ and $\delta''$.  Right: The intermediate region, bounded by $|y-1|=\delta$ with $y=1+\frac t {\sqrt n}$, and lying outside of the neighborhoods of the turning points $t=\pm 2\sqrt a$. The gray parts are where the intermediate asymptotics and the turning point asymptotics match. }
\label{fig:domains}
\end{figure}

Indeed, taking  $\delta$ in Lemma \ref{prop-approx-at-0} to be $1-\delta'$ with sufficiently small $\delta'$, we see that \eqref{approx-not-at-1} holds   for $y\in D_1: |y|<1-\delta'$. On the other hand, from
\eqref{approx-out-Charlier-formula}, putting in use Stirling's formula, we have  \eqref{approx-not-at-1}
for $y\in D_2: |y-[0, 1]| >\delta''$. With appropriately small $\delta'$ and $\delta''$, we have $D_1\cup D_2$ covers  $|y-1|\geq \delta$; see Figure \ref{fig:domains} for an illustration. For $y$ in a neighborhood of the positive real axis, one may have to use the identity $(-1)^n \frac{\Gamma(n-x)}{\Gamma(-x)}=\frac {\Gamma(x+1)} {\Gamma(1+x-n)}$.
It is worth mentioning that $\varepsilon_1$   is also exponentially small for finite $x$; cf. the derivation leading to  Corollary \ref{cor:zeros}.

It is asked in Bo and Wong \cite{bo-wong1994} whether there exists a uniform asymptotic expansion for $C_n^{(a)}(n\beta)$ in the interval $-\infty< \beta \leq \delta$, where the constant $\delta\in (0, 1)$. Theorem \ref{thm-approx-not-at-1} provides   an answer.

\subsection{An intermediate  region }\label{subsec:intermediate}
Now we have explicitly derived an asymptotic formula of  $C^{(a)}_n(ny)$ for $|y-1|\geq \delta$. The domain of uniformity contains a neighborhood  of infinity, and a neighborhood of $y=0$, that is, an end-point of the equilibrium measure.
What left is an $O(1)$ neighborhood of $y=1$, namely, $|y-1|< \delta$.
We will see that $y=1$ is of significance since  a pair of turning points coalesce there.

The last two decades see dramatic changes in difference equation methods in this respect. For example, Wang and Wong \cite{wang-wong2003,wang-wong2005} established a turning point theory for second order linear difference equations, the theory is further completed by several authors, including Cao and Li \cite{cao-li2014}; see also the review article \cite{wong2014}. In the mentioned works, special functions, such as the Airy functions and Bessel functions, are employed to describe the asymptotic behavior  at the turning point. However, there is a connection problem to be solved, namely, one has  to further determine  which solution behaves as this  asymptotic solution. In many cases, solving such a  connection problem turns out to be  a hard question, in particular  using the difference equation methods alone. An attempt  to overcome the difficulty has been made by Geronimo \cite{geronimo2009}. Some of the material here also appears in Huang \cite{huang2018}.

In the present paper, we demonstrate, using the Charlier polynomials as an example, how to solve the connection problem. We use three types of asymptotics, respectively in  non-oscillatory region, intermediate  region (with constants to be determined), and at the turning point (in the forms of a linear  combination), all obtained via difference equation methods. Matching adjacent regions outside in, we determine asymptotics in the inner regions, step by step.

Substituting
$C_n^{(a)}(x)=(2a)^{\frac n 2} \frac {\Gamma\left((n+1)/ 2\right )} {\Gamma\left( 1/ 2\right )}P_n(x)$ into \eqref{charlier-difference-eq} gives the symmetric canonical form
\begin{equation}\label{charlier-difference-eq-symmetric}
P_{n+1}(x)-\left (A_n x+B_n\right ) P_n(x)+P_{n-1}(x)=0,
\end{equation}
where, as $n\to\infty$,
\begin{equation}\label{difference-eq-balance}
A_n=\frac 1 {\sqrt{2a}}\frac {\Gamma\left (\frac {n+1} 2\right )} {\Gamma\left(\frac {n} 2+1\right )} \sim \frac 1 {\sqrt n} \sum_{s=0}^\infty \frac {\alpha_s} {n^s},~~  B_n= - \frac {n+a} {\sqrt{2a}} \frac {\Gamma\left(\frac {n+1} 2\right )} {\Gamma\left(\frac {n} 2+1\right )}\sim \sqrt n \sum_{n=0}^\infty \frac {\beta_s}{n^s},
\end{equation}with $\alpha_0=\frac 1{\sqrt a}$, $\alpha_1=-\frac 1 {4\sqrt a}$,  $\beta_0= -\frac 1 {\sqrt a}$, $\beta_1= \frac 1 {4\sqrt a} -\sqrt a$,  and $\beta_s=-\alpha_s-a \alpha_{s-1}$ for $s=1,2,\cdots$.

Now we introduce a new local variable $t$ at $y\sim 1$,
  \begin{equation}\label{transf-at-1-Charlier}
   x=n\left ( 1+\frac t {\sqrt n}\right ).
  \end{equation}
In what follows, we  consider the case when $t$ is taken away from the real interval $[-\sqrt n, 2\sqrt a]$, to which all    zeros belong; cf.  Remark \ref{rem:not-the-same} for the upper bound $t=2\sqrt a$.
An  idea   in \cite{huang-cao-wang2017} applies, with modifications, so that in such an outer region, \eqref{charlier-difference-eq-symmetric} possesses a pair of non-vanishing asymptotic solutions  of
  the form
 \begin{equation}\label{approx-out-Charlier}
 P_n(x) \sim  \exp\left ( \sqrt n  \phi_{-1}(t) + \phi_0(t) +  \sum^\infty_{k=1} \frac {\phi_k(t)}{n^{k/2}}\right ),
 \end{equation}where $\phi_k(t)$, to be determined,  are functions independent of $n$.

 \begin{rem}
In earlier literature, such as Wong and Li \cite[Eq.(1.5)]{wong-li-1992}, there is an extra factor of the form $n^\alpha$ attached to the expression in \eqref{approx-out-Charlier}.
In the present case, we have $x=n\left(1+\frac t{\sqrt n}\right )$, and
we can always write
\begin{equation*}
n^\alpha= x^\alpha \exp\left( -\alpha\log\left(1+  \frac t{\sqrt n}\right )\right )= x^\alpha \exp\left( \sum^\infty_{k=1}\frac {\alpha(-1)^k t^k} k \frac 1 {n^{k/2}}\right )ㄛ
\end{equation*}
in which $x^\alpha$ behaves like a constant, and the new exponential function combines with the one in \eqref{approx-out-Charlier}.
That explains  why we can ignore the possible factor $n^\alpha$.
\end{rem}

Here, unlike \cite{huang-cao-wang2017}, instead of   working out  an infinite series,  we focus on the leading terms $\phi_{-1}(t)$ and $\phi_0(t)$.
From \eqref{transf-at-1-Charlier},  we see  that $t$ shifts for the same $x$ when $n$ varies,
and we write
\begin{equation}\label{x-invariant}
x=n\left ( 1+\frac t {\sqrt n}\right )=(n+1)\left ( 1+\frac {t_+} {\sqrt {n+1}}\right )=(n-1)\left ( 1+\frac {t_-} {\sqrt {n-1}}\right ).
\end{equation}The shifts $t_\pm-t$ can be expanded in descending powers of $n$; see Section \ref{sec:Intermediate} for details.
In view of \eqref{x-invariant},  substituting \eqref{approx-out-Charlier} into \eqref{charlier-difference-eq-symmetric}, and
equalizing the constant terms   on both sides, we have the first order differential equation
 \begin{equation}\label{differential-eq-phi--1-Charlier}
e^{- \phi_{-1}'(t)}+e^{ \phi_{-1}'(t)}=\frac t {\sqrt a}
\end{equation}
 for $\phi_{-1}$, which implies  $e^{ \phi_{-1}'(t)}=\frac {t\pm \sqrt{t^2-4a}}{2\sqrt a}$, where the branch is chosen so that $\sqrt{t^2-4a}$ is analytic in $\mathbb{C}\setminus [-2\sqrt a, 2\sqrt a\;]$ and being real positive for $t>2\sqrt a$. We pick the minus sign first, and obtain a solution to \eqref{differential-eq-phi--1-Charlier},
\begin{equation*}
\phi_{-1}(t)=t\log\frac {t-\sqrt{t^2-4a}}{2\sqrt a} +\sqrt{t^2-4a}+C_{-1},
\end{equation*}such that $e^{ \phi_{-1}'(t)}=\frac {t-\sqrt{t^2-4a}}{2\sqrt a}$,
where $
C_{-1}$ is a constant independent of $n$. Substituting the expression into \eqref{charlier-difference-eq-symmetric} and further equalizing the coefficients of $1/\sqrt n$ gives
\begin{equation*}
\phi_0'(t)= \frac a {\sqrt{t^2-4a}}+\frac 1 2 \sqrt{t^2-4a}  - \frac t{ 2(t^2-4a)}   +\frac 1 2 C_{-1}.
\end{equation*}
Thus we have
\begin{equation*}
\phi_0(t)=-\frac 1 4  \log(t^2-4a)   +\frac 1 4 t \sqrt{t^2-4a}  +\frac 1 2 C_{-1}t +C_0,
\end{equation*}where $
C_0$ is a constant.
With such  $\phi_{-1}(t)$ and $\phi_0(t)$, the asymptotic solution \eqref{approx-out-Charlier} now reads
\begin{equation}\label{approx-out-Charlier-at-1}
P_n(x)  \sim C \exp\left ( \sqrt n  \left [ t\log\frac {t-\sqrt{t^2-4a}}{2\sqrt a} +\sqrt{t^2-4a} \right ]   -\frac 1 4  \log(t^2-4a)  +\frac t 4 \sqrt{t^2-4a}   \right ),
\end{equation}where   $x=n  ( 1+\frac t {\sqrt n}  )$, and
$C=C(n, t)=
e^{ \sqrt n   C_{-1}  +\frac 1 2 C_{-1}t +C_0}$.

Similarly, if we take the alternative choice $e^{ \phi_{-1}'(t)}=\frac {t+ \sqrt{t^2-4a}}{2\sqrt a}$, we have the other formal (asymptotic) solution
\begin{equation}\label{approx-out-Charlier-at-1-2nd}
\tilde P_n(x)  \sim \tilde C \exp\left ( -\sqrt n  \left [ t\log\frac {t-\sqrt{t^2-4a}}{2\sqrt a} +\sqrt{t^2-4a} \right ]   -\frac 1 4  \log(t^2-4a) -\frac t 4 \sqrt{t^2-4a}   \right ),
\end{equation}where $\tilde C=\tilde C(n, t)=
e^{ -\sqrt n   \tilde C_{-1}  -\frac 1 2 \tilde C_{-1}t -\tilde C_0}$, with $\tilde C_{-1}$ and $ \tilde C_0$ being constants.

Now we match the asymptotic approximations \eqref{approx-out-Charlier-formula} and
  \begin{equation}\label{Charlier-at-1-inter}
  C_n^{(a)}(x)\sim  (2a)^{\frac n 2} \frac {\Gamma\left((n+1) /2\right )} {\Gamma\left(1/ 2\right )}\left ( A(x) P_n(x)+B(x)\tilde P_n(x)\right )~~\mbox{as}~~n\to\infty,
\end{equation}
with $x=ny=n(1+t/\sqrt n)$,  at the transition area described as
\begin{equation*}
n^{1/6}\ll |t|\ll  n^{1/4}.
\end{equation*}
Here $A(x)$ and $B(x)$, to be determined,  depend only on $x$ and $a$, but not on $n$.

First, we see that for $t\in (2\sqrt a, +\infty)$, $\tilde P_n(x)$ is dominant and $P_n(x)$ is recessive. The approximation in  \eqref{approx-out-Charlier-formula} is also  recessive, hence the coefficient $B(x)$ vanishes.
The remaining coefficient in \eqref{Charlier-at-1-inter} and the constants $C_{-1}$ and $C_0$ can be determined by the matching process to give $A(x)=\left (\frac {\Gamma(x+1)}{a^x}\right )^{1/2}$,  $C_{-1}=0$ and  $C_0=-\frac 3 4\log 2-\frac 1 4\log\pi +\frac a 2$. As a result, we have
\begin{thm}
\label{thm-approx-intermediate}For $x=ny=n(1+\frac t {\sqrt n})$, with $|y-1|<\delta$ and $|t-(-\infty, 2\sqrt a\,]|> r$, $\delta$ and $r$ being positive constants, it holds
\begin{equation}\label{approx-intermediate}
C_n^{(a)}(x)\sim \frac C{\sqrt{w(x)}}  \exp\left ( \sqrt n  \left [ t\log\frac {t-\sqrt{t^2-4a}}{2\sqrt a} +\sqrt{t^2-4a} \right ]   -\frac 1 4  \log(t^2-4a)  +\frac t 4 \sqrt{t^2-4a}   \right )
\end{equation}as $n\to\infty$, where $w(x)= \frac {a^x}{\Gamma(x+1)}$ and  $C=(2a)^{\frac n 2} \frac {\Gamma\left( (n+1)/ 2\right )} {\Gamma\left(1/ 2\right )} 2^{-\frac  3  4} \pi^{- \frac 1 4}
e^{ \frac a  2}$.

\end{thm}
\begin{rem}
\label{rem:not-the-same}
 In Huang-Cao-Wang \cite{huang-cao-wang2017}, an assumption for the symmetric difference equation \eqref{charlier-difference-eq-symmetric} is
$B_n\sim \sum \frac {\beta_s}{n^s}$. Here in the Charlier case we have $B_n\sim \sqrt n \sum \frac {\beta_s}{n^{s}}$  instead; cf. \eqref{difference-eq-balance}.  Nevertheless, after the change of variable \eqref{transf-at-1-Charlier}, which is nonlinear in $n$, we can still apply the method in \cite{huang-cao-wang2017} to extract the leading terms. In this intermediate case the method seems to be superior to the one illustrated in the previous subsection.

It is worth pointing out that  altering the assumption has significant impact to the turning point analysis that follows: The Charlier case seems to go beyond the framework of the turning point theory; see Wong \cite{wong2014}.  The general assumption in \cite{wong2014} is the difference equation \eqref{charlier-difference-eq-symmetric} with
\begin{equation*}
A_n\sim n^{-\theta}\sum^\infty_{s=0}\frac {\alpha_s}{n^s} ,~~B_n\sim \sum^\infty_{s=0} \frac {\beta_s}{n^s},
\end{equation*}just as in \cite{huang-cao-wang2017}. The types of asymptotic solutions are described  by the characteristic  equation $\lambda^2-(\alpha_0 y+\beta_0)\lambda+1=0$, $x=n^{\theta} y$, which has two roots that coincide when  $\alpha_0 y+\beta_0=\pm 2$. Such a $y$   is a turning point; see \cite{wong2014}. However,  in the Charlier case, the $P_n$ term in \eqref{charlier-difference-eq-symmetric} is dominant, and the characteristic equation degenerates to $\alpha_0 y+\beta_0=0$ with $y=n^{-1}x$, giving rise to a single critical point at $y=1$.
Yet in local variable $t$, introduced in \eqref{transf-at-1-Charlier}, the characteristic  equation takes the form $\lambda^2-\frac t {\sqrt a} \lambda+1=0$, which locates two turning points at $t=\pm 2\sqrt a$.
\end{rem}

\subsection{At turning points}

In the terminology of
  Baik et al.\;\cite{baik-et-al2007}, there is a saturated-band-void configuration in a shrinking neighborhood of $y=1$. Here saturated region is defined as   an open subinterval of   maximal
length in which  the equilibrium measure realizes the upper constraint.  A void by definition  is an open subinterval of  maximal length in which  the equilibrium measure realizes the lower constraint, namely, $0$. The band lies in between.

As mentioned in Remark \ref{rem:not-the-same}, for the Charlier difference equation \eqref{charlier-difference-eq-symmetric}, the index equation degenerates. With the local re-scaling  \eqref{transf-at-1-Charlier}, we are capable of determining the band in variable $t$ as $t\in (-2\sqrt a, 2\sqrt a)$. We proceed to derive the transition asymptotics of  $C_n^{(a)}(x)$  in the band, or, more challenging, at the turning points $t=\pm 2\sqrt a$, where $x=n(1+\frac t {\sqrt n})$.

Bo and Wong \cite{bo-wong1994} used Bessel functions to give the approximation of  $C^{(a)}_n(ny)$
for $0<\varepsilon \leq y\leq M<\infty$, restricted to the real axis.   In \cite{ou-wong2010}, Ou and Wong addressed global asymptotic formulas using Riemann-Hilbert approach.
The uniform expansions they derived involve a quite complicated combination of the Airy functions. Tedious calculation is needed if one has to draw local leading asymptotics
from their  formulas.

In what follows, we apply the theory of Wang and Wong (see \cite{wong2014}) to obtain the asymptotic form of solutions at turning points.   The most attention will be paid to connect the asymptotic solution with the Charlier polynomials,  by matching  the approximation  here with those in earlier subsections.

Still, we start from the standard form difference equation \eqref{charlier-difference-eq-symmetric}, relaxing the initial conditions.  As mentioned in Remark \ref{rem:not-the-same}, the fundamental assumptions on $A_n$ and $B_n$
are not fulfilled, as compared with \cite{wang-wong2003, wong2014}. However, by introducing the local transformation
$x=n(1+\frac t {\sqrt n})$, the method of Wang and Wong keeps applicable.

As a first step, we treat the turning point $t=t_0=2\sqrt a$. The idea in \cite{wang-wong2003} is to seek an asymptotic solution to \eqref{charlier-difference-eq-symmetric}, of the form
\begin{equation*}
Q_n(x)=\sum^\infty_{s=0}  \chi_s(\xi) \delta^s;
\end{equation*}as suggested by Costin and Costin \cite{costin-costin1996}, where $\delta=n^{-\alpha}$,  $\xi=n^\sigma \eta(t)$, with $\eta(t_0)=0$, $\eta(t)$ being a one-to-one mapping in a neighborhood of $t=t_0$. Substituting the expressions in the difference equation and matching the leading terms, we determine $\sigma=\frac 1 3$ and $\alpha=\frac 1 2$. Ignoring  lower order terms in $n$, we have
\begin{equation*}
 \chi_0(\xi_+)-(A_n x+B_n)\chi_0(\xi)+\chi_0(\xi_-)=0,
\end{equation*}
where $\xi_+=(n+1)^{1/3} \eta(t_+)$ and $ \xi_-=(n-1)^{1/3} \eta(t_-)$ vary around $\xi$,  $t_\pm$ being given  in \eqref{x-invariant}.
Repeatedly using Taylor expansions, constantly ignoring lower order terms, and taking advantage of the symmetry form, we further obtain the differential equation
\begin{equation*}
 \chi_0''(\xi)=c^3 \xi \;  \chi_0(\xi),
\end{equation*}where the constant $c=a^{-1/6}/\eta'(t_0)$. Similarly, other $\chi_k$ solve inhomogeneous equations of the same type as $\chi_0$. This indicates the involvement of Airy functions to represent all $\chi_k$.  As in \cite{wang-wong2003}, we may write the asymptotic solution in the more accurate form
\begin{equation}\label{airy-type-sol}
Q_n(x)=\chi\left (n^{1/3}\eta +n^{-1/6}\Phi\right ) \sum^\infty_{k=0}  \frac {A_k(\eta)} {n^{k/2}}   + n^{-1/6} \chi'\left (n^{1/3}\eta +n^{-1/6}\Phi\right ) \sum^\infty_{k=1}  \frac {B_k(\eta)} {n^{k/2}},
\end{equation}and proceed to determine $\eta(t)$, $\Phi(t)$ and $A_0(\eta)$.
The function $\chi(\xi)$ in \eqref{airy-type-sol}  solves the Airy equation.    Here
we have used the fact that  $\eta(t)$ defines a conformal mapping at $t=t_0$, therefore, the non-oscillating coefficients $A_k$ and $B_k$, as functions of $t$, can be regarded as functions of $\eta$.
Substituting \eqref{airy-type-sol} into \eqref{charlier-difference-eq-symmetric}, with details given in Section \ref{sec:Airy-type-approx}, we do have
$\eta(t)$, a conformal mapping at $t=2\sqrt a$, being positive
  for $t>2\sqrt a$, such that
\begin{equation}\label{conformal-mapping-introduction}
\frac 2 3 \left ( \eta(t) \right )^{3/2} =t\log\frac {t+\sqrt{t^2-4a}}{2\sqrt a}-\sqrt{t^2-4a},~~t\in \mathbb{C}\setminus (-\infty, 2\sqrt a].
\end{equation}We can further determine
\begin{equation}\label{Phi-sol}
\sqrt {\eta(t)} \Phi(t)=-\frac {t\sqrt{t^2-4a}} 4,
\end{equation}where $\sqrt {\eta(t)}$ and    $\sqrt{t^2-4a}$ are positive for $t> 2\sqrt a$.
Also, $A_0(\eta)$ is solved up to a constant factor independent of both $n$ and $t$. We pick
\begin{equation}\label{A-0-sol}
A_0(\eta)= \left ( \frac {t^2-4a} {4a\eta}\right )^{-1/4}.
\end{equation}
Choosing  the Airy function $\chi$ in \eqref{airy-type-sol} to be $\Ai$ and $\Bi$, and denoting the asymptotic solution respectively as $Q_n(x)$ and $\tilde Q_n(x)$, we have
\begin{equation*}
C_n^{(a)}(x)\sim (2a)^{\frac n 2} \frac {\Gamma\left((n+1)/ 2\right )} {\Gamma\left( 1/ 2\right )}\left (   K_1(x) Q_n(x)+K_2(x) \tilde Q_n(x)\right );
\end{equation*}cf.\;\eqref{turn-point-right}.
The coefficients  $K_1(x)$ and $K_2(x)$ are determined by matching the approximation with the intermediate asymptotic behavior \eqref{approx-intermediate}, as conducted in Section \ref{subsec:turning-point-right}.
\begin{thm}\label{thm-right-turning-point}
For $x=n\left (1+\frac t {\sqrt n}\right )$, in  an $O(1)$  neighborhood of the turning point $t=2 \sqrt a$, as $n\to\infty$, it holds
\begin{eqnarray}\label{approximation-turn-point-right}
 C_n^{(a)}(x)&=&  C_{K, n} x^{ \frac 1 {12}} \left (\frac {a^x}{\Gamma(x+1)}\right )^{-\frac 1 2} \left ( \frac {t^2-4a} {4a\eta(t)}\right )^{-1/4} \left\{ \Ai \left (n^{\frac 1 3}\eta(t) +n^{-\frac 1 6}\Phi(t)\right )\left [1+O\left (n^{-\frac 1 2}\right )\right ]\right . \nonumber \\
   & &\left . +\Ai' \left (n^{\frac 1 3}\eta(t) +n^{-\frac 1 6}\Phi(t)\right )O\left (n^{-\frac 2 3}\right )\right\},
\end{eqnarray}
where $\eta(t)$ and $\Phi(t)$ are given in \eqref{conformal-mapping-introduction} and \eqref{Phi-sol}, and  $C_{K, n}=(2a)^{\frac n 2} \frac {\Gamma\left((n+1)/ 2\right )} {\Gamma\left( 1/ 2\right )} \left ( \frac \pi{2a}\right )^{1/4}e^{\frac a 2}$.
\end{thm}

Now we turn to the other turning point $t=-2\sqrt a$.
In this case,
we substitute
 \begin{equation*}
 C_n^{(a)}(x)=(-1)^n (2a)^{\frac n 2} \frac {\Gamma\left((n+1)/ 2\right )} {\Gamma\left( 1/ 2\right )}\mathcal{Q}_n(x)
 \end{equation*}
 into \eqref{charlier-difference-eq}, to  give the canonical difference equation
\begin{equation}\label{charlier-difference-eq-symmetric-minus-introduction}
\mathcal{Q}_{n+1}(x)+\left (A_n x+B_n\right ) \mathcal{Q}_n(x)+\mathcal{Q}_{n-1}(x)=0.
\end{equation}The factor $(-1)^n$ is brought in for a technical reason: The coefficient
$A_nx+B_n=t/\sqrt a+O(1/\sqrt n)$ in \eqref{charlier-difference-eq-symmetric-minus-introduction} at $t=-2\sqrt a$, corresponds to $-(A_nx+B_n)=-t/\sqrt a+O(1/\sqrt n)$ in \eqref{charlier-difference-eq-symmetric}  at $t=2\sqrt a$, so that the derivation leading to Theorem \ref{thm-right-turning-point} may be used, with minor modifications, in this case.

We still seek asymptotic solutions to \eqref{charlier-difference-eq-symmetric-minus-introduction}
of the form \eqref{airy-type-sol}, with $\tilde\eta$, $\tilde\Phi$, $\mathcal{A}_k$ and $\mathcal{B}_k$ taking the places of $\eta$, $\Phi$, $A_k$ and $B_k$, respectively. Solving differential equations yields
\begin{equation}\label{conformal-mapping-minus-introduction}
\frac 2 3 {\tilde\eta}^{3/2} =t\log \frac{-t+\sqrt{t^2-4a}}{2\sqrt a}+\sqrt{t^2-4a},
\end{equation}such that ${\tilde\eta}(t)$ is a conformal mapping at $t=-2\sqrt a$ with ${\tilde\eta}(-2\sqrt a)=0$, ${\tilde\eta}(t)>0$ for $t<-2\sqrt a$, and  ${\tilde\eta}(t)\sim   a^{-1/6} (-2\sqrt a-t )$  for $t\sim -2\sqrt a$. Here, as before,   $\sqrt{t^2-4a}$ is analytic in $\mathbb{C}\setminus [-2\sqrt a, 2\sqrt a]$, and behaves like $t$ as $t\to \infty$. The logarithm takes principal branch.

Also, one obtains
\begin{equation}\label{Phi-sol-left}
 \sqrt {\tilde\eta} {\tilde\Phi} =\frac  1 4 t \sqrt{t^2-4a},
\end{equation}where ${\tilde\Phi}(t)$ is analytic in a neighborhood of $t=-2\sqrt a$, and $\sqrt {{\tilde\eta}(t)}$ take positive values for $t<-2\sqrt a$.
The leading coefficient $\mathcal{A}_0({\tilde\eta})$ can be determined up to a factor independent of $n$ and  $t$. We pick
\begin{equation}\label{A-0-sol-left}
\mathcal{A}_0({\tilde\eta})= \left ( \frac {t^2-4a} {4a{\tilde\eta}}\right )^{-1/4}= \exp\left (\frac {\pi i} 2-\frac 1 4 \log  (t^2-4a) +\frac 1 4\log (4a{\tilde\eta})\right ),
\end{equation}so that $\mathcal{A}_0({\tilde\eta})$ is analytic in a neighborhood of $t=-2\sqrt a$, being real positive for $t<-2\sqrt a$. Here in the neighborhood, $t$ is an analytic function of ${\tilde\eta}$.

Now we determine the asymptotic approximation in a upper half neighborhood of $t=-2\sqrt a$, again by a matching process.
 Hence the solution to \eqref{charlier-difference-eq-symmetric-minus-introduction}, corresponding to the Charlier polynomials, has the asymptotic approximation $\mathcal K_1(x) \mathcal Q_n(x)+\mathcal K_2(x) \tilde {\mathcal Q}_n(x)$; cf. \eqref{turn-point-left}, where
\begin{equation*}
\mathcal{Q}_n(x)\sim {\mathcal{A}_0({\tilde\eta})} \Ai \left (\omega n^{\frac 1 3}\left ({\tilde\eta}  +\frac {{\tilde\Phi}}{\sqrt n}\right )\right ) ,~~\tilde {\mathcal{Q}}_n(x)\sim {\mathcal{A}_0({\tilde\eta})} \Ai \left (\omega^2  n^{\frac 1 3}\left ({\tilde\eta}  +\frac {{\tilde\Phi} }{\sqrt n}\right )\right ),
\end{equation*}with $\omega =e^{ 2\pi i/3 }$, and with
  $\mathcal K_1(x)$ and $\mathcal K_2(x)$ to be determined by matching \eqref{turn-point-left} with the intermediate asymptotic formula \eqref{approx-intermediate}.  The matching is carried out first in the sector $\arg \tilde\eta\in (-\pi, -\pi/3)$, namely $\arg (t+2\sqrt a)\in (0, 2\pi/3)$, where $\tilde {\mathcal Q}_n(x)$ is dominant over $\mathcal{Q}_n(x)$; cf. Figure \ref{fig:domains} for an illustration of the sector,  and  see Section \ref{subsec:turning-point-left} for a full description of the process. Accordingly, we have

\begin{thm}\label{thm-left-turning-point}
For $x=n\left (1+\frac t {\sqrt n}\right )$, in  a small    neighborhood of the turning point $t=-2 \sqrt a$, as $n\to\infty$, it holds
 \begin{equation}\label{approximation-turn-point-left-real}
C_n^{(a)}(x)\sim \frac {C_{\mathcal{K}, n} x^{\frac 1 {12}}\mathcal{A}_0({\tilde\eta})}{\sqrt{w(x)} }
  \left [ \cos(x\pi) \Ai \left (  n^{\frac 1 3}\left ({\tilde\eta}  +\frac {{\tilde\Phi}}{\sqrt n}\right )\right )- \sin(x\pi)  \Bi \left ( n^{\frac 1 3}\left ({\tilde\eta}  +\frac {{\tilde\Phi}}{\sqrt n}\right )\right )\right ] ,
\end{equation}
where $w(x)=\frac {a^x}{\Gamma(x+1)}$, ${\tilde\eta}(t)$ and ${\tilde\Phi}(t)$ are given in \eqref{conformal-mapping-minus-introduction}  and \eqref{Phi-sol-left}, $\mathcal{A}_0({\tilde\eta})$ takes the branch as in \eqref{A-0-sol-left}, and  $C_{\mathcal{K}, n}= (-1)^n (2a)^{\frac n 2} \frac {\Gamma\left((n+1)/ 2\right )} {\Gamma\left( 1/ 2\right )}  \left (\frac \pi {2a}\right )^{1/4}e^{a/2}$.
\end{thm}
\vskip .5cm

It is worth mentioning that $e^{\mp {\pi i} /3}\Ai \left ( s e^{\pm 2\pi i/3}\right )=\frac 1 2\left ( \Ai(s)\mp i \Bi(s)\right )$; cf. \cite[Eq.(9.2.11)]{nist}. Thus the functions  in the square brackets can be rewritten as
\begin{equation*}
 e^{-\left(x\pi+  \pi /3\right )i} \Ai \left (\omega\Theta  \right )+e^{\left (x\pi+  \pi /3\right )i} \Ai \left (\omega^{-1}\Theta \right ) ,
\end{equation*}where $\Theta=n^{\frac 1 3}\left ({\tilde\eta(t)}  +  {{\tilde\Phi(t)}}/{\sqrt n}\right )$. For $t+2\sqrt a$ small with $\Im (t+2\sqrt a)>0$, we have
 \begin{equation*}
 \left |e^{-\left(x\pi+\pi /3\right )i} \Ai \left (\omega\Theta  \right )
\right |= O\left (n^{-1/12}\right ) e^{\sqrt n \left [\Im (t+2\sqrt a) +O\left( (t+2\sqrt a)^{3/2}\right )\right ]}
 \end{equation*}
 and
 \begin{equation*}
 \left |e^{ \left(x\pi+\pi /3\right )i} \Ai \left (\omega^{-1}\Theta  \right )
\right |=O\left (n^{-1/12}\right )
  e^{\sqrt n \left [-\Im (t+2\sqrt a) +O\left( (t+2\sqrt a)^{3/2}\right )\right ]}.
 \end{equation*}
Here  use has been made of the  asymptotic approximation
$\Ai(z)\sim \frac 1 {2\sqrt\pi} z^{-\frac 1 4} e^{-\frac 2 3z^{3/2}}$ as $|\arg z|<\pi$, $z\to\infty$, and the facts that  ${\tilde\eta}^{3/2}\sim a^{-1/4} (-2\sqrt a-t )^{3/2}$ and    $\tilde\Phi(t)\sim a^{5/6}$ for $ t\sim -2\sqrt a$, as drawn from  \eqref{conformal-mapping-minus-introduction} and \eqref{Phi-sol-left}.
 Therefore, it holds
 \begin{equation*}
 \left |e^{-\left(x\pi+\frac \pi 3\right )i} \Ai \left (\omega\Theta  \right )
\right |\gg
 \left |e^{ \left(x\pi+\frac \pi 3\right )i} \Ai \left (\omega^{-1}\Theta  \right )
\right |
 \end{equation*}
for $\Im (t+2\sqrt a)>0$  and  $ t+2\sqrt a$ small such that $\Im (t+2\sqrt a)>\left|O\left( (t+2\sqrt a)^{3/2}\right )\right |$. Picking up the main contribution of \eqref{approximation-turn-point-left-real} for such $t$, one has
\begin{equation}\label{approximation-turn-point-left}
C_n^{(a)}(x)\sim C_{\mathcal{K}, n} \frac { x^{\frac 1 {12}}\mathcal{A}_0({\tilde\eta})}{\sqrt{w(x)} }  e^{-(x\pi+  \pi /3)i} \Ai \left (\omega n^{\frac 1 3}\left ({\tilde\eta}  +\frac {{\tilde\Phi}}{\sqrt n}\right )\right ) .
\end{equation}Similarly, in the small neighborhood, with $\Im(t+2\sqrt a)<0$, instead of $e^{-\left(x\pi+\frac \pi 3\right )i} \Ai \left (\omega\Theta  \right )$, one retains  $e^{ \left(x\pi+  \pi/ 3\right )i} \Ai \left (\omega^{-1}\Theta  \right )$.

The asymptotic approximation \eqref{approximation-turn-point-left} has an error as   in \eqref{approximation-turn-point-right},
 of the form $O(n^{-2/3})\Ai'$, which has no impact on the leading   asymptotic location of zeros.
 Also, for \eqref{approximation-turn-point-left-real} to become an equality, a factor $(1+O(n^{1/2}))$
 should be attached to both $\Ai$ and $\Bi$.

Later in Section \ref{subsec:turning-point-left}, a coherence check is made. Approximations for $t<-2\sqrt a$ and $t>-2\sqrt a$ are derived and are shown in consistence with   \eqref{approx-in-Charlier-formula} and \eqref{approx-between-via-intermediate}, obtained respectively from the non-oscillatory region asymptotics, and intermediate asymptotics.

\begin{rem}
The asymptotic forms \eqref{approximation-turn-point-right} and \eqref{approximation-turn-point-left} differ from \cite{wong2014} in the appearance of a shift $n^{-1/6}\Phi(t)$ in the variable, yet they do agree with \cite{wang-wong2003}. For sure one can re-expand  the Airy function by using \eqref{chi-Taylor} and \eqref{chi-prime-Taylor} to get rid of $\Phi$.  However, by this way the leading term can be extracted appropriately; see the absence of $B_0$ in \eqref{airy-type-sol}.


\end{rem}

\subsection{Discussion and arrangement of the rest of the paper}
The Charlier polynomials
$C_n^{(a)}(n y)$ display an interesting asymptotic
behavior near this interval $y\in [0, 1]$,  which is the support of the relevant  uniform density-$1$  equilibrium
measure, where
the zeros of the polynomial asymptotically accumulate.  On the   interval,  the zeros of $C_n^{(a)}(n y)$ are in
one-to-one correspondence with, and for large $n$ are extremely close to, the corresponding nodes of
orthogonality, namely the values $y = k/n$, $k = 0, 1, 2, \cdots$.  However, the region $y > 1$ is outside the support
of the equilibrium measure, so there are asymptotically no zeros. The density of zeros jumps suddenly from
$1$ to $0$ near $y = 1$.

This is an unusual situation, because normally there is an interval of continuous transition between the
upper and lower constraints for the equilibrium measure.
The reason, as demonstrated in the present paper, is the existence of a shrinking   band region,
described by   an appropriate
re-scaling near $y = 1$, namely $y=1+t/\sqrt n$, with the transition region $-2\sqrt a < t<2\sqrt a$  in terms of the variable
$t$,  the endpoints of which are simple
turning points for the three-term recurrence relation.
     Thus the band asymptotics, involving Airy functions and stated  in Theorems \ref{thm-right-turning-point} and \ref{thm-left-turning-point}, furnishes a dramatic yet smooth  transition.
   Near the other endpoint  $y=0$, there is another transition described by a simple  ratio of Gamma functions, as given in Theorem \ref{thm-approx-not-at-1}.

The motivation of the present investigation is twofold. First, from a difference equation point of view, the Charlier polynomials are special, as indicated in   the basic  assumptions \eqref{difference-eq-balance}, in that they can not fit into any of the known cases; cf. Wong \cite{wong2014}, see also Remark \ref{rem:not-the-same}.
Therefore, it is desirable to derive uniform and non-uniform asymptotics of the polynomials via difference equations. As described in the present section, we have put in use various methods, all of a difference equation nature, to deal with asymptotic approximations respectively in the outside region, an intermediate region, and near the turning points. The overlapping   domains of validity  actually cover the whole complex $y$-plane. Here $y=x/n$ is  the re-scaling that normalizes the support of the equilibrium measure to $y\in [0,1]$. In particular, we obtain uniform asymptotic expansions at a pair of coalescing turning points in the $y$-plane.

Another motivation is that the Charlier polynomials may serve as a  model in the study of the   Heun class equations.  For example, when one considers the connection problems between fundamental solutions for the confluent Heun's equation (CHE), and the doubly-confluent Heun's equation (DHE),  a central piece seems to be a three term recurrence relation satisfied by the coefficients of a certain  Frobenius solution, essentially  similar to \eqref{charlier-difference-eq}, in that they share the same equilibrium support $y\in [0,1]$, after a re-scaling   $y=x/n^2$. Here $x$ is the accessory parameter in CHE or DHE. Most likely, the techniques used here, such as those in Section \ref{subsec:origin}, could play a role. Eigenvalue problems and   root polynomials for CHE and DHE also seem to be relevant.

It is worth mentioning that in a paper of  Dai-Ismail-Wang \cite{dai-ismail-wang}, non-oscillatory asymptotic approximations are derived, and matched with asymptotic behavior resulted from the turning point techniques of Wang and Wong (see \cite{wong2014}). The present investigation repeatedly makes use of matching processes, just as \cite{dai-ismail-wang} did, yet the treatment of  coalescing turning points (in Section \ref{sec:Airy-type-approx}) and the uniform approximation at the end-point $y=0$ (in Section \ref{subsec:origin}) seem to be novel, and applicable elsewhere.  It is also worth noticing that in Section \ref{subsec:turning-point-right}-\ref{subsec:turning-point-left}, the matching processes are closely related to the Stokes lines of the asymptotic solutions; see also Remark \ref{rem:subdominant-sector}.

The rest of the paper is arranged as follows.  In Section \ref{subsec:non-oscillatory-asymptotics}, we derive a uniform asymptotic approximation for
  $C_n^{(a)}(n y)$ as $|y- [0,1]|>r$ and $n\to\infty$, $r$ being a generic  positive constant. In Section \ref{subsec:origin}, we provide uniform asymptotics
for
  $C_n^{(a)}(n y)$ as $|y|<r$ and $n\to\infty$, where $r\in (0, 1)$. A combination of the results covers the $y$-domain $|y-1|>\delta$, and thus proves   Theorem \ref{thm-approx-not-at-1}.
Then,  Section   \ref{sec:Intermediate} is devoted to the asymptotic approximation of $C_n^{(a)}(x)$ in the intermediate region described as $|t\pm 2\sqrt a|> r$ and $|y-1|<\delta$, where $x=ny=n\left ( 1+\frac t{\sqrt n}\right)$, and again $r$ and $\delta$ are generic positive constants. The result turns out to be Theorem \ref{thm-approx-intermediate}. In  Section \ref{sec:Airy-type-approx}, we determine two uniform asymptotic approximations at the turning
points $t=2\sqrt a$ and $t=-2\sqrt a$, respectively, and we prove Theorems \ref{thm-right-turning-point} and \ref{thm-left-turning-point} in this section.
 In the last section, Section \ref{sec:comparison}, we apply Theorems \ref{thm-approx-not-at-1}, \ref{thm-right-turning-point} and \ref{thm-left-turning-point}, to obtain local behaviors near the end point $y=0$, and turning points  $t=2\sqrt a$ and $t=-2\sqrt a$, respectively. The results are compared with the asymptotic formulas obtained earlier
 by Bo and Wong \cite{bo-wong1994} and Goh \cite{goh1998}.

\section{Non-oscillatory regions, the origin, and proof of Theorem \ref{thm-approx-not-at-1} }

\subsection{Non-oscillatory asymptotics}\label{subsec:non-oscillatory-asymptotics}
We are in a position to prove \eqref{approx-out-Charlier-formula}. To this aim, we need to show the validity of \eqref{w-k-Charlier} beforehand.
It is appropriate to write  \eqref{w-k-Charlier} as
\begin{equation}\label{w-k-Charlier-error}
w_k(x):=(x-k) \left \{ 1+  \frac {1- a}  {x-k} - \frac {ak} {(x-k)^2}   +\varepsilon_k \right \}:=(x-k) \left \{ 1+\delta_k\right \}.
\end{equation}
We   estimate the error terms as follows.

\begin{lem}\label{lemma-w-k-estimate} Assume that $x=n y$ with   $|y- [0,1]|>r$ for a certain positive constant $r$. Then, there exist positive constants $M_0$, $M_1$ and $N$, such that
\begin{equation}\label{delta-k-Charlier}
\left | \delta_k \right | \leq \frac {M_0  } {n }
\end{equation} and
\begin{equation}\label{epsilon-k-Charlier}
\left | \varepsilon_k \right | \leq \frac {M_1} {n^2}
\end{equation}for $k=1,2,\cdots, n$ and $n>N$.
 \end{lem}

 \noindent {\bf Proof}: We prove the lemma by induction in $k$. Comparing  $w_1(x)=x-a$ with \eqref{w-k-Charlier-error},
 we see that
 \begin{equation*}
 \left |\delta_1\right |\leq
\frac  {|1-a|/r} n~~\mbox{and}~~  \left |\varepsilon_1\right |\leq \frac {a/r^2}{n^2},
 \end{equation*} as long as $|y-[0,1]|>r$.

 Assume the validity of \eqref{delta-k-Charlier} for   $k-1$,  we proceed to  show that  \eqref{epsilon-k-Charlier} holds for index $k$.  Indeed, substituting \eqref{w-k-Charlier-error} into \eqref{w-k-Charlier-first-order} with $a_k=k+a$ and $b_k=ak$, we have
\begin{equation*}
 (x-k)\varepsilon_k=\left (\frac {ak} {x-k}-\frac {a(k-1)}{x-k+1}\right ) +\frac {a(k-1)}{x-k+1}\left(1- \frac 1 {1+\delta_{k-1}}\right ).
\end{equation*}Straightforward estimation then gives
\begin{equation*}
\left | \varepsilon_k\right |\leq \frac{ a\left ( 1+r+2rM_0\right )} {r^3} \frac 1 {n^2}
\end{equation*}for $k\leq n$ and $n>N=2M_0$. Accordingly, from \eqref{w-k-Charlier-error} we have
\begin{equation*}
\left |\delta_k\right |\leq  \left (\frac {|1-a|} r+\frac {3 a}{r^2} +\frac a {r^3}\right ) \frac 1 n
\end{equation*}for $k\leq n$ and $n>N=2M_0$.
Hence, assigning
\begin{equation*}
 M_0= \frac {|1-a|} r+\frac {3 a}{r^2} +\frac a {r^3},~~     M_1= \frac{ a\left ( 1+r+2rM_0\right )} {r^3}~~\mbox{and}~~  N=2M_0,
\end{equation*}we complete the proof of the lemma.\qed
\vskip .3cm

Now substituting \eqref{w-k-Charlier} into \eqref{approx-out-Charlier-monic}, we obtain  the asymptotic approximation
\eqref{approx-out-Charlier-formula}  for $y$ keeping a distance $r$ from $[0,1]$. Here use has been  made   of the trapezoidal rule to give
\begin{equation*}
\sum^n_{k=1}\log\left (y-\frac{k}{n} \right ) = n\int^1_0\log(y-t) dt +\frac 1 2\log\frac{y-1} y+O\left (\frac 1 n\right ),
\end{equation*}and the right-hand side terms explicitly give
\begin{equation*}
n\left [  y \log  y -(  y-1)\log( y-1)-1\right ]+\frac 1 2\log\frac{y-1} y.
\end{equation*}
We also pick up the later terms in \eqref{w-k-Charlier}  by approximating
\begin{equation*}
\sum_{k=1}^n \log \left \{ 1+\frac {(1-a)x- k }{(x-k)^2}\right \}\sim \sum_{k=1}^n  \frac {(1-a)x- k }{(x-k)^2},
\end{equation*}which in turn can be approximated by
\begin{equation*}
\int^1_0  \frac {(1-a)  y-t }{(  y-t)^2}dt = \log\frac y {y-1}-\frac {a}{y-1},
\end{equation*}each time with an error $O(1/n)$.

As mentioned in Szeg\H{o} \cite[p.395]{szego-book},   forming the real part of the approximation of a polynomial $P_n(x)$, the asymptotic formula for $\frac 1 2 P_n(x)$ arises.
Following Szeg\H{o},
a
formal derivation from \eqref{approx-out-Charlier-formula}   gives a formula for fixed $y\in (0, 1)$, namely,
\begin{eqnarray}\label{approx-in-Charlier-formula}
  C_n^{(a)}(n y)    &=& 2 n^{n} \sqrt{\frac y {1-y }}  \exp\left (-\frac {a}{y-1}\right ) \exp\left\{n  \left [   y  \log\frac {  y }{ 1- y } -1\right ] \right\} (y-1)^n \nonumber\\
    & & \times \left \{ -\sin\left (n y\pi \right )\left [1+O\left(\frac 1 n\right )\right ]+  \cos\left (n y\pi \right ) O\left(\frac 1 n\right )\right \}.
\end{eqnarray}
A similar idea  has been  applied to, e.g., \cite{huang-cao-wang2017}.  It is readily verified that \eqref{approx-in-Charlier-formula} agrees with Theorem \ref{thm-approx-not-at-1} for $y\in (0, 1)$.

\subsection{Asymptotics for $C^{(a)}_n(x)$ as $|x|<n$}\label{subsec:origin}

More precisely, we derive a uniform asymptotic formula for $C^{(a)}_n(x)$ as $\frac {|x|} n\leq r<1$ and $n\to\infty$, $r$ being an arbitrary  constant. To do so,
we   apply  the method of dominant balance to the three-term  recurrence formula \eqref{charlier-difference-eq}; cf. Bender and Orszag \cite[Ch.5]{bender-orszag-1999-book} for the method.
Indeed,  by matching the  $ C^{(a)}_n$ term with $C^{(a)}_{n-1}$ we obtain
  an asymptotic solution of \eqref{charlier-difference-eq},
 \begin{equation*}
 C_s(n)=C_s(n,x) \sim (-a)^n  \frac {\Gamma(n+1)}{\Gamma(n-x+1) }  \exp\left( -\frac{an}{n-x}\right ) ,
 \end{equation*} for large $n$ and large $n-x$.
While   a comparison of the    $C^{(a)}_{n+1}$ and $C^{(a)}_n$ terms gives rise to another
 asymptotic solution
\begin{equation*}
C_d(n)=C_d(n,x) \sim (-1)^n \Gamma (n-x).
 \end{equation*} These are the only possible balances.
 It is readily seen that $C_d(n,x)$ is dominant over  $C_s(n,x)$ as $n\to\infty$  for finite $x$.

Refinements are available for both solutions. Indeed, we have
\begin{lem}\label{lemma:uniformity}
There exist solutions $C_d(n,x)$ and $C_s(n,x)$ of \eqref{charlier-difference-eq}, such that
\begin{equation}\label{Q-n-asym-homogeneous}
 C_d(n,x) = (-1)^n  \Gamma(n-x)\left [ \sum^{N-1}_{k=0} \frac {\varphi_k(y)} {n^k}+O\left(\frac 1 {n^N}\right )\right ],~~~n\to\infty,
 \end{equation}and
 \begin{equation}\label{C-s-def}
C_s(n,x) = (-a)^n  \frac {\Gamma(n+1)}{\Gamma(n-x+1) }  \left [ \sum^{N-1}_{k=0} \frac {\psi_k(y)} {n^k}+O\left(\frac 1 {n^N}\right )\right ],~~~n\to\infty,
 \end{equation}
 where all functions $\varphi_k(y)$ and $\psi_k(y)$ are independent of  $n$ and analytic in the unit $y$-disc, and the error terms $O\left(  1 / {n^N}\right )$, $N=1,2,\cdots$, are uniform
   for   $|y|\leq r$, with $x=ny$ and arbitrary $r\in (0, 1)$.
\end{lem}

\noindent {\bf Proof}:
We prove \eqref{Q-n-asym-homogeneous} in full detail. First, taking the places of $C^{(a)}_n(x)$ in \eqref{charlier-difference-eq} with
$(-1)^n  \Gamma(n-x)w_n(x)$, we obtain  an equation
\begin{equation}\label{w-difference-eq}
w_{n+1}(x)-\left(1+\frac a {1-y}\frac 1 n\right )w_n(x)
+\left(\sum^\infty_{k=1} \frac a {(1-y)^{k+1}}\frac 1 {n^k}\right ) w_{n-1}(x)=0,
\end{equation}where $x=ny$. Justifying \eqref{Q-n-asym-homogeneous} amounts to show that \eqref{w-difference-eq} has a uniform approximate solution
\begin{equation}\label{w-n-expansion}
 w_n(x)=\sum^{N-1}_{k=0} \frac {\varphi_k(y)} {n^k}+O\left(\frac 1 {n^N}\right ):=W_N(n, x)+\varepsilon_N(n, x).
\end{equation}

 We  will see that the leading coefficient $\varphi_0(y)$ is determined up to a constant factor.
To this end, we write $x=ny=(n+1) y_+=(n- 1) y_-$, and we see that the shifts
  \begin{equation}\label{shift-y}
   y_\pm- y=\mp\frac y n \left (1\pm \frac 1 n\right )^{-1}=
y\sum^\infty_{k=1} \frac {(\mp 1)^k} {n^k},
\end{equation}for finite $y$ and large $n$.
Substituting the formal solution \eqref{w-n-expansion} into \eqref{w-difference-eq}, expanding $\varphi_k(y_+)$ and $\varphi_k(y_-)$ into Taylor expansions at $y$, and picking up the $O(1/n)$ terms, we have
\begin{equation*}
  - y     \varphi_0'(y)
 +\frac { ay} {(1-y)^2}  \varphi_0(y) = 0.
\end{equation*}Solving this first order differential equation gives an analytic function in the unit disc
\begin{equation*}
\varphi_0(y)=  \exp \left (  \frac {ay}{1-y}\right ),
\end{equation*}
up to a constant factor.

The functions $\varphi_k(y)$ for $k=1,2,\cdots$   can be  determined   recursively  by comparing  the $O(1/n^{k+1})$ terms. Assume that we have had
analytic functions $\varphi_k(y)$  for $k=0,1,\cdots, N-1$, we proceed to show that $\varphi_N(y)$ is uniquely determined by the analyticity.

Indeed, substituting \eqref{w-n-expansion} into \eqref{w-difference-eq}, we obtain
\begin{equation}\label{DE-for-reminder}
\varepsilon_N(n+1, x)-\left(1+\frac a {1-y}\frac 1 n\right )\varepsilon_N(n, x)
+\left(\sum^\infty_{k=1} \frac a {(1-y)^{k+1}}\frac 1 {n^k}\right )\varepsilon_N(n-1, x)=R_N(n,x),
\end{equation}
 where
 \begin{equation*}
                R_N(n,x)=
                -W_N({n+1},x)+\left(1+\frac a {1-y}\frac 1 n\right )W_N({n },x )
-\left(\sum^\infty_{k=1} \frac a {(1-y)^{k+1}}\frac 1 {n^k}\right )W_N({n-1},x)
\end{equation*}takes the form
\begin{equation}\label{R-n-expansion}
R_N(n,x)  =  \sum^\infty_{k=N+1}\frac{r_k(y)}{n^k},
\end{equation}with each function $r_k(y)$, independent of $n$, being analytic in $|y|<1$, and satisfying $|r_k(y)|\leq M \rho^k$ as $|y|\leq r$ for arbitrary $r\in (0, 1)$, and    generic positive constants $M$ and $\rho$. The expansion \eqref{R-n-expansion}
needs justification. First, the absence of all $n^{-k}$ terms for $k=0,1,\cdots, N-1$ is clear since that is how those $\varphi_k$'s are determined. The analyticity and bounds for $r_k(y)$ can be achieved by straightforward calculation. A typical example is
\begin{equation*}
\varphi_0(y_+)-\varphi_0(y)
=\sum^\infty_{k=1}\frac{\varphi_0^{(k)}(y)}{k!}(y_+-y)^k
=\sum^\infty_{s=1}\left \{\sum^s_{k=1}\frac{y^k \varphi_0^{(k)}(y)}{k!}\frac{k(k+1)\cdots (s-1)}{(s-k)!}\right\}\frac {(-1)^s}{n^s}.
\end{equation*}Here use has been made of the first equality of \eqref{shift-y}, and the binomial formula. It is readily seen that the coefficients of $n^{-s}$ share the same domain of analyticity as  $\varphi_0(y)$.
What is more, for a constant $r_1\in (r, 1)$, one has a bound
$|\varphi_0(y)|\leq M$ for some (generic) $M$ as $|y|\leq r_1$, which implies
\begin{equation*}
\frac{ \left  | \varphi_0^{(k)}(y)\right |}{k!}\leq \frac {Mr_1}{(r_1-r)^{k+1}}~~\mbox{for}~~|y|<r
\end{equation*}by Cauchy's integral formula. In turn we have a bound of the form $M\rho^s$ for the coefficient of $n^{-s}$, that is
\begin{equation*}
\left |\sum^s_{k=1}\frac{y^k \varphi_0^{(k)}(y)}{k!}\frac{k\cdots (s-1)}{(s-k)!}\right |\leq \frac{M r_1}{r_1-r}\sum^s_{k=1}\left (\frac{r}{r_1-r}\right )^{k} \frac{k\cdots (s-1)}{(s-k)!}=\frac {M r }{ r_1-r } \left(\frac {r_1} {r_1-r}\right )^{s}.
\end{equation*}

Bearing in mind that $\varepsilon_N(n, x)$ has a leading term $\varphi_N(y)/n^N$, in view of \eqref{R-n-expansion}, and picking up all $1/n^{N+1}$ terms in \eqref{DE-for-reminder}, we obtain an equation
\begin{equation}\label{fisrt-order-eq}
-y\varphi_N'(y)+\left ( \frac{ay}{(1-y)^2}-N\right )\varphi_N(y)=r_{N+1}(y),
\end{equation}to which all solutions are
\begin{equation*}
\varphi_N(y)=\left [C-\int^y_0 \eta^{N-1} e^{-\frac a {1-\eta}} r_{N+1} (\eta) d\eta\right  ] y^{-N} e^{\frac a {1-y}}.
\end{equation*}For $N\geq 1$, there is only one solution, namely when the constant $C=0$, that is analytic at the origin, and in the unit disc as well,
owing to the analyticity of $r_{N+1}(y)$ for $|y|<1$.

Now we turn back to \eqref{DE-for-reminder}, to give a rigorous proof for the existence and uniformity of the asymptotic solution. More precisely, we shall show that for appropriate $N$, \eqref{DE-for-reminder} has a solution satisfying
\begin{equation*}
\left|\varepsilon_N(n, y)\right |\leq M/n^N~~~\mbox{for}~~|y|\leq r,~~n\to\infty,
\end{equation*}where $r\in (0, 1)$, and again $M$ is a generic constant.
To this end,
we follow several steps   in
   Wong and Li \cite[Sec.\;3]{wong-li-1992}. The first step is to write
 \eqref{DE-for-reminder} as the following
 contractive form
 \begin{equation}\label{DE-contractive}
\varepsilon_N(n+1, x)-\varepsilon_N(n, x)=R_N(n,x)+  \frac a {n-x} \varepsilon_N(n, x)
- \frac {an} {(n-x)(n-x-1)} \varepsilon_N(n-1, x).
\end{equation}
One can formally derive from \eqref{DE-contractive} the equation
\begin{equation}\label{DE-sol-contra}
 \varepsilon_N(n, x)=\sum^\infty_{k=n} \left [-R_N(k,x)-  \frac a {k-x} \varepsilon_N(k, x)
+ \frac {ak} {(k-x)(k-x-1)} \varepsilon_N(k-1, x)\right ].
\end{equation}It is readily verified that every solution of \eqref{DE-sol-contra} is a solution of \eqref{DE-contractive}.

As in \cite{wong-li-1992},    we solve \eqref{DE-sol-contra}
  by using  the method of  successive approximation. Define the sequence $\{h_s(n,x)\}$ by $h_0(n,x)=0$ and
 \begin{equation*}
  h_{s+1}(n,x)= \sum^\infty_{k=n} \left [-R_N(k,x)-  \frac a {k-x} h_{s}(k,x)
+ \frac {ak} {(k-x)(k-x-1)} h_{s}(k-1,x)\right ]
 \end{equation*}for $n\geq n_0$, and  $h_{s+1}(n,x)\equiv 0$ for $n<n_0$.

From \eqref{R-n-expansion}  we have   $|R_N(n,x)|\leq M /n^{N+1}$ for $n\geq n_0$, with $n_0$ big enough, and $M$ may depend on $N$ but not on $n$ and $y$ for $|y|<r$. It is worth noting that   $N$ is fixed, but  could be large.  Hence for  $n\geq n_0$  with $n_0>N$, we have
\begin{equation*}
 \left | h_{1}(n,x)\right |=\left | \sum^\infty_{k=n}  R_N(k,x) \right |\leq
M \sum^\infty_{k=n}  \frac 1 {k^{N+1}}\leq \frac {2M} N\frac 1 {n^N}.
 \end{equation*}
 For $n\geq n_0$, $|x|<rn_0$ and $n_0$ large enough, it also holds
\begin{equation*}
 \left | h_2(n,x)-h_{1}(n,x)\right |\leq   \frac {2a M}{(1-r) N}\sum^\infty_{k=n}\frac 1 {k^{N+1}}+  \frac {2a M}{(1-r)^2 N}\sum^\infty_{k=n}\frac 1 {(k-1)^{N+1}}  \leq
  \frac {2LM} {N^2 n^N},
 \end{equation*}where the constant $L=\frac {2a} {1-r}+\frac {4a}{(1-r)^2}$, and we have used the fact  that    $\frac 1 {(n-1)^N}<\frac 2 {n^N}$ for $n\geq n_0$.  The process can be repeated. By induction, for  $n\geq n_0$ and  $|x|<rn_0$, it holds
 \begin{equation*}
 \left | h_{s+1}(n,x)-h_{s}(n,x)\right |\leq
  \frac {2L^sM} {N^{s+1} n^N},~~~s=1,2,\cdots.
 \end{equation*}
 Therefore, if $N$ is chosen such that $N>L$, we see that
 \begin{equation*}
 \varepsilon_N(n, x)=\sum_{s=0}^\infty \left \{ h_{s+1}(n,x)-h_{s}(n,x)\right \}
 \end{equation*}is absolutely and uniformly convergent for $n\geq n_0$ and $|x|<r n_0$, such that
 \begin{equation*}
 \left |\varepsilon_N(n, x)\right |\leq \frac {2M}{N-L} \frac 1 {  n^N}.
 \end{equation*}
Thus it furnishes a solution to
 \eqref{DE-sol-contra}, and hence solves \eqref{DE-contractive}. We have proved \eqref{Q-n-asym-homogeneous}.

 The derivation of  \eqref{C-s-def} is entirely parallel, and the asymptotic results are not used in the present paper, so  we skip the proof here.
\qed
\vskip 1cm

Now that    $C_d(n,x)$ and  $C_s(n,x)$ are  linear independent    solutions of \eqref{charlier-difference-eq} in the unit $y$-disc with $x=ny$, there exist functions $K_d(x)$ and  $K_s(x)$, depending only on $x$, but not on $n$, such that
\begin{equation}\label{disc-connection}
 C^{(a)}_n(x)= K_d(x) C_d(n,x)+ K_s(x) C_s(n,x),~~~x=ny,~~~|y|<1.
 \end{equation}
 For $|y|\leq r<1$, with both $n$ and $n-x=n(1-y)$ being large,
we rewrite  \eqref{approx-out-Charlier-formula}   as
$
C^{(a)}_n(x)  \sim  (-1)^n e^{ \frac {an }{n-x }}\frac {\Gamma\left (n-x\right )}{\Gamma(-x)}$ by applying Stirling's formula,
and  compare it with \eqref{Q-n-asym-homogeneous}. Then we see a perfect matching
\begin{equation}\label{C-n-combination}
C^{(a)}_n(x)\sim \frac{e^a}{\Gamma(-x)} C_d(n,x)
\end{equation}for $y$ in compact subsets of the unit disc cut along  the segment $[0, 1]$, with $x=ny$.
Here, by perfect matching we mean that the later terms on both sides
   of \eqref{C-n-combination},   determined via the same recurrence formula \eqref{charlier-difference-eq}, are exactly the same.
That is, the error term
 \begin{equation*}
 \varepsilon_1=K_s(x) C_s(n,x),
 \end{equation*}
 is beyond all orders, in other words, is exponentially small as compared with the right-hand side term in \eqref{C-n-combination}.
 Therefore, we have

\begin{lem}\label{prop-approx-at-0} There is a large-$n$ asymptotic approximation
\begin{equation}\label{Charlier-at-0}
C^{(a)}_n(x)  = (-1)^n \frac {\Gamma\left (n-x\right )}{\Gamma(-x)} \left [ \sum^{N-1}_{k=0} \frac {\varphi_k(y)} {n^k}+O\left(\frac 1 {n^N}\right )\right ]+ \varepsilon_1,
\end{equation} as $|y| <\delta$ for  arbitrary constant $\delta\in (0, 1)$,  where $x=ny$,  $\varphi_0(y)=e^{ \frac {a  }{1-y }}$, and
each $\varphi_k(y)$, $k=1,2,\cdots$, is the unique analytic solution of \eqref{fisrt-order-eq}    in $|y|<1$.
The error terms  $O\left(  1 /{n^N}\right )$ for  $N=1, 2,\cdots$ is uniform in $|y|<\delta$, and
$\varepsilon_1=K_s(x) C_s(n,x)$ is exponentially small   as compared with the leading right-hand side term in the domain $|y|<\delta$ and  $|y-[0,1] |> r $, for arbitrary positive $r$.
 \end{lem}

A sharper bound of $\varepsilon_1$ can be obtained for finite $x$. Indeed, for $|x|<R$ with arbitrary $R$ and $n$ large enough, straightforward calculation gives
\begin{equation*}
\left | \frac {K_s(x) C_s(n,x)}{C_d(n,x)}\right |\leq M e^{- n \log n+(1+\log a) n  +( 2R+1/2)\log n  },
\end{equation*}where $M$ is
a constant independent of $n$. Here the function $K_s(x)$, no matter how it is determined, is bounded for bounded $x$,  as follows from the boundedness of all other functions in \eqref{disc-connection}, and the non-vanishing of $C_s(n,x)$. This implies that for   $|x|<R$, the error $\varepsilon_1$ in
\eqref{Charlier-at-0} is exponentially small as compared with the leading contribution. Hence the smallest zeros are determined by $  1 /{\Gamma(-x)}$, with exponentially small perturbations. Then, as a corollary of Lemma \ref{prop-approx-at-0}, we reproduce the following well-known result (see, e.g., \cite[Thm.3]{ou-wong2010}).

\begin{cor}\label{cor:zeros}
Denote the zeros of $C^{(a)}_n(x)$ in increasing order as
$x_{1,n}<x_{2,n}<\cdots<x_{n,n}$. Then, for any fixed $k\in \mathbb{N}=\{1,2,3,\cdots\}$, it holds
 \begin{equation}\label{zeros}
 x_{k,n}=k-1+\mbox{an~exponentially~small~term},~~\mbox{as}~n\to\infty.
 \end{equation}
 \end{cor}

\section{ Intermediate asymptotics   and proof of Theorem \ref{thm-approx-intermediate}}
\label{sec:Intermediate}

The  difference equation \eqref{charlier-difference-eq-symmetric} has asymptotic solutions of the form \eqref{approx-out-Charlier} in the intermediate region, namely, in $|y-1|<\delta$ and $|t\pm 2\sqrt a| >\tilde \delta$, where $x=ny=n\left (1+\frac t {\sqrt n}\right )$, $\delta$ and $\tilde \delta$ are constants. Now we take the leading terms. More precisely, we determine $\phi_{-1}$ and $\phi_0$, so that
$P_n(x)\sim e^{\sqrt n \phi_{-1}(t) +\phi_0(t)}$. The derivation is similar to  \cite{huang-cao-wang2017}, with modifications.

Temporarily, we assume that $\left |t-(-\infty,  2\sqrt a\; ]\right | >\tilde \delta$.
From \eqref{transf-at-1-Charlier} and \eqref{x-invariant}, for the same $x$  when $n$ varies,  we write
\begin{equation*}
x=n\left ( 1+\frac t {\sqrt n}\right )=(n+1)\left ( 1+\frac {t_+} {\sqrt {n+1}}\right )=(n-1)\left ( 1+\frac {t_-} {\sqrt {n-1}}\right ),
\end{equation*}and see that, for large $n$,
\begin{equation*}
t_+=-\frac 1 {\sqrt{n+1}}+\sqrt{\frac n {n+1}} t=t-\frac 1 {\sqrt n}-\frac t {2n} +O\left ( n^{-3/2}\right )+O\left (\frac t{n^2}\right ),
\end{equation*}and
\begin{equation*}
t_-= \frac 1 {\sqrt{n-1}}+\sqrt{\frac n {n-1}} t=t+\frac 1 {\sqrt n}+\frac t {2n} +O\left ( n^{-3/2}\right )+O\left (\frac t{n^2}\right ).
\end{equation*}
Substituting \eqref{approx-out-Charlier} into \eqref{charlier-difference-eq-symmetric} gives
\begin{equation*}
\begin{array}{l}
\displaystyle{e^{\sqrt{n+1}\phi_{-1}(t_+) +\phi_0(t_+) -\sqrt{n }\phi_{-1}(t) -\phi_0(t) +O(1/n)}+e^{\sqrt{n-1}\phi_{-1}(t_-) +\phi_0(t_-) -\sqrt{n }\phi_{-1}(t) -\phi_0(t) +O(1/n)}} \\
    =\alpha_0 t+\frac {\alpha_1+\beta_1}{\sqrt n}+O(1/n).
\end{array}
\end{equation*}
Expanding $\phi_{-1}$ and $\phi_0$ at $t$, we have
\begin{equation*}
\begin{array}{l}
\displaystyle{e^{-\phi'_{-1}(t)+\frac 1 {2\sqrt n} \left [ \phi_{-1}(t)-t \phi'_{-1}(t)+\phi''_{-1}(t)-2\phi'_0(t)\right ]  +O  (\frac 1 n  )}+e^{\phi'_{-1}(t)+\frac 1 {2\sqrt n} \left [ -\phi_{-1}(t)+t \phi'_{-1}(t)+\phi''_{-1}(t)+2\phi'_0(t)\right ]  +O  (\frac 1 n  )}}
 \\
     =\frac t{\sqrt a}-\frac {\sqrt a}{\sqrt n}+O\left (\frac 1 n\right ).
\end{array}
\end{equation*}
Equalizing the coefficients of $1$ and $1/\sqrt n$, respectively,  on both sides, we have the first order differential equations \eqref{differential-eq-phi--1-Charlier}, namely,
$e^{- \phi_{-1}'(t)}+e^{ \phi_{-1}'(t)}=\frac t {\sqrt a}$,
 for $\phi_{-1}$,  and
\begin{equation}\label{differential-eq-phi-0-Charlier}
       \left (  e^{\phi_{-1}' } -  e^{-\phi_{-1}' }\right )   \phi_0' = -\sqrt a - e^{-\phi_{-1}' }\left [ \frac  1 2\phi_{-1} -\frac {t} 2 \phi_{-1}'  + \frac 1 2\phi_{-1}''  \right ] -
       e^{\phi_{-1}' }\left [ -\frac  1 2\phi_{-1} +\frac {t} 2 \phi_{-1}'  + \frac 1 2\phi_{-1}''  \right ]
\end{equation}for $\phi_0$.

From  \eqref{differential-eq-phi--1-Charlier} one sees that $e^{ \phi_{-1}'(t)}=\frac {t\pm \sqrt{t^2-4a}}{2\sqrt a}$, each corresponds to an asymptotic solution of \eqref{charlier-difference-eq-symmetric}. For example, taking  the minus sign and integrating, we have
\begin{equation*}
\phi_{-1}(t)=t\log\frac {t-\sqrt{t^2-4a}}{2\sqrt a} +\sqrt{t^2-4a}+C_{-1},
\end{equation*}such that $e^{ \phi_{-1}'(t)}=\frac {t-\sqrt{t^2-4a}}{2\sqrt a}$,
where $
C_{-1}$ is a constant independent of $n$. Substituting the expression into \eqref{differential-eq-phi-0-Charlier}, we obtain
\begin{equation*}
\phi_0'(t)= \frac a {\sqrt{t^2-4a}}+\frac 1 2 \sqrt{t^2-4a}  - \frac t{ 2(t^2-4a)}   +\frac 1 2 C_{-1}.
\end{equation*}
Then the asymptotic solution \eqref{approx-out-Charlier-at-1} is readily derived.

Now we match the asymptotic approximations \eqref{approx-out-Charlier-formula} and \eqref{Charlier-at-1-inter}  with $x=ny$ and $y=1+t/\sqrt n$,  at the transition area described as
\begin{equation*}
n^{1/6}\ll |t|\ll  n^{1/4}.
\end{equation*}
As mentioned in Section \ref{subsec:intermediate}, we see from \eqref{approx-out-Charlier-formula} or \eqref{approx-not-at-1} that $C_n^{(a)}(n+\sqrt n t )$ is recessive for $t\in (2\sqrt a, +\infty)$.
The asymptotic solution $\tilde P_n(x)$ in \eqref{approx-out-Charlier-at-1-2nd} is dominant as compared with $P_n(x)$. Hence the coefficient $B(x)$ vanishes. Thus $C_n^{(a)}(x)\sim (2a)^{\frac n 2} \frac {\Gamma\left((n+1) /2\right )} {\Gamma\left(1/ 2\right )}A(x)P_n(x)$.

Substituting in $y=1+\frac t {\sqrt n}$ with $|t|\ll  n^{1/4}$, we can write \eqref{approx-out-Charlier-formula} as
\begin{equation*}
C^{(a)}_n\left (x \right )\sim  \left (\frac n e\right )^n n^{\frac 1 4} t^{-\frac 1 2} \exp\left (-\frac {a\sqrt n} t-\sqrt n \, t\log t+\frac t 2 \sqrt n \log n+\sqrt n \, t+\frac {t^2} 2-\frac {t^3}{6\sqrt n}\right ).
\end{equation*}
On the other hand, a combination of \eqref{approx-out-Charlier-at-1} with \eqref{Charlier-at-1-inter} as $|t| \gg n^{1/6}$ gives
\begin{equation*}
C_n^{(a)}(x)\sim A(x) C_K(n, t)  \exp\left (  \sqrt n  \left [ t\log\sqrt a-t\log t+t -\frac a t  \right ]   -\frac 1 2  \log t  +\frac {t^2} 4 -\frac  a 2   \right ),
\end{equation*}where   $x=n  ( 1+\frac t {\sqrt n}  )$, and
$C_K(n, t)= (2a)^{\frac n 2} \frac {\Gamma\left((n+1) /2\right )} {\Gamma\left(1/ 2\right )} e^{ \sqrt n   C_{-1}  +\frac 1 2 C_{-1}t +C_0}$.
A comparison of these formulas, with the aid   of Stirling's formula, gives
\begin{equation*}
A(x)\sim 2^{-\frac 1 2 } e^{-C_{-1} \sqrt n-\frac 1 2 C_{-1} t-C_0+\frac a 2} x^{\frac 1 4} e^{\frac 1 2 x\log x-\frac x 2} a^{-\frac x 2}
\end{equation*}for large $n$. Here use has been made of the fact that
$\frac 1 2 x\log x-\frac x 2\sim \frac 1 2 n\log n-\frac n 2 +\frac t 2 \sqrt n\log n+\frac {t^2 } 4$ and $x^{1/4}\sim n^{1/4}$. We also  see the  involvement  of $\phi_1(t)$; see \eqref{approx-out-Charlier}. The function   $A(x)$ is independent of $n$, so we have $C_{-1}=0$, and we may take
\begin{equation*}
A(x)=(w(x))^{-\frac 1 2} =\left (\frac {\Gamma(x+1)}{a^x}\right )^{1/2}~~\mbox{and}~~C_0=-\frac 3 4\log 2-\frac 1 4\log\pi +\frac a 2,
\end{equation*}where $w(x)=\frac {a^x}{\Gamma(x+1)}$ is the discrete weight shown in \eqref{Charlier-orthogonal}. Thus completes the proof of \eqref{approx-intermediate} in Theorem \ref{thm-approx-intermediate}.

The same idea leading to \eqref{approx-in-Charlier-formula} applies. From  \eqref{approx-intermediate} we have
\begin{equation}\label{approx-between-via-intermediate}
C_n^{(a)}(x)\sim \frac {2C}{\sqrt{w(x)}\; (4a-t^2)^{1/4}}  \cos\left (2\sqrt {a   n } \left [\sin\theta -\theta\cos\theta \right ]+a\sin\theta\cos\theta -\frac \pi 4    \right ),
\end{equation}where $x=n(1+\frac t {\sqrt n})$ with $t=2\sqrt a\cos\theta \in (-2\sqrt a, 2\sqrt a)$, $\theta\in (0, \pi)$, and as in \eqref{approx-intermediate}, $w(x)= \frac {a^x}{\Gamma(x+1)}$ and  $C=(2a)^{\frac n 2} \frac {\Gamma\left(  (n+1)/2\right )} {\Gamma\left(  1 /2\right )} 2^{- \frac 3 4} \pi^{- \frac 1 4}
e^{\frac  a  2}$.

\section{Airy-type approximations }\label{sec:Airy-type-approx}
Adapting   the change of variable
$x=n(1+\frac t {\sqrt n})$,
we consider  the
  turning points  $t=t_0=\pm 2\sqrt a$.
The aim is to determine the asymptotic solution \eqref{airy-type-sol}, or, the leading term
\begin{equation*}
Q_n(x)\sim {A_0(\eta)}\;  \chi\left (n^{1/3}\eta(t) +n^{-1/6}\Phi(t)\right ) ,
\end{equation*}namely, to determine $\eta(t)$, $\Phi(t)$ and $A_0(\eta)$. Here $\chi$ solves the Airy equation $\chi''(\tau)=\tau \chi(\tau)$.
The derivation is very similar to Wang and Wong \cite{wang-wong2002}. But, since the basic assumption \eqref{difference-eq-balance}, and the form in \eqref{airy-type-sol}, differ from those in  \cite{wang-wong2002},  we need to describe the key steps briefly.

We begin with $Q_{n\pm 1}(x)$, preserving $x$ while $n$ varies. We write $x=n+ {\sqrt n}\; t=(n\pm 1)+ {\sqrt {n\pm 1}}\; {t_\pm}$
as in \eqref{x-invariant}. For large $n$, we see that the shifts
\begin{equation*}
t_+-t\sim -\frac 1 {\sqrt n}-\frac t {2n}~~\mbox{and}~~t_--t\sim  \frac 1 {\sqrt n}+\frac t {2n}.
\end{equation*}
Also we turn to the variable of $\chi$. To this aim, we denote
\begin{equation}\label{chi-variable-shift}
(n\pm 1)^{1/3}\eta(t_\pm) +(n\pm1)^{-1/6}\Phi(t_\pm):=n^{1/3}\eta(t) +n^{-1/6}\Phi(t)+ n^{-1/6}u_\pm.
\end{equation}It is readily seen that both $u_\pm$ admit asymptotic expansion of the form $\sum^\infty_{k=0} \frac {u_k}{n^{k/2}}$, with the leading terms given as
\begin{equation}\label{u+-leading}
u_+\sim \sum^\infty_{k=0} \frac {u_{+,k}}{n^{k/2}}=-\eta'(t) +\left (\frac {\eta(t)} 3-\frac {t\eta'(t)} 2+\frac {\eta''(t)} 2-\Phi'(t)\right )\frac 1 {\sqrt n}+\cdots
\end{equation}
and
\begin{equation}\label{u--leading}
u_-\sim \sum^\infty_{k=0} \frac {u_{-,k}}{n^{k/2}}= \eta'(t) +\left (-\frac {\eta(t)} 3+\frac {t\eta'(t)} 2+\frac {\eta''(t)} 2+\Phi'(t)\right )\frac 1 {\sqrt n}+\cdots.
\end{equation}

 Regarding $\zeta$ and  $\mu$ as   free variables, we take a closer look at $\chi(n^{1/3} \zeta +n^{-1/6} \mu)$.
\begin{lem}\label{lem:chi}
Assume  that $\chi$ solves the Airy equation. Then the equality holds
\begin{equation}\label{chi-Taylor}
\chi\left (n^{1/3} \zeta +n^{-1/6} \mu\right )= \chi\left(n^{1/3} \zeta \right )X(n; \zeta, \mu)+ n^{-1/6}   \chi'\left(n^{1/3} \zeta\right )Y(n; \zeta, \mu),
\end{equation}
where the coefficients possess the asympototic expansions
\begin{equation}\label{X-Y-series}
 X(n; \zeta, \mu)=\sum^\infty_{k=0} \frac {X_k(\zeta, \mu)} {n^{k/2}}~~\mbox{and}~~
Y(n; \zeta, \mu)=\sum^\infty_{k=0} \frac {Y_k(\zeta, \mu)} {n^{k/2}},
 \end{equation}with
\begin{equation}\label{X-0-Y-0-sol}
  X_0(\zeta, \mu)=\frac 1 2\left (e^{\sqrt \zeta\mu}+ e^{-\sqrt \zeta\mu}\right ),~~~~ Y_0(\zeta, \mu)=\frac 1 {2\sqrt\zeta}\left (e^{\sqrt \zeta\mu}- e^{-\sqrt \zeta\mu}\right ),
 \end{equation}
 \begin{equation}\label{X-k-sol}
 X_k(\zeta, \mu)=\frac 1 {2\sqrt\zeta} \int^\mu_0 s X_{k-1}(\zeta, s)\left ( e^{\sqrt\zeta (\mu-s)}- e^{\sqrt\zeta (s-\mu)} \right ) ds,~~k=1,2,\cdots
 \end{equation}
  and
   \begin{equation}\label{Y-k-sol}
Y_k(\zeta, \mu)=\frac 1 {2\sqrt\zeta} \int^\mu_0 s Y_{k-1}(\zeta, s)\left ( e^{\sqrt\zeta (\mu-s)}- e^{\sqrt\zeta (s-\mu)} \right ) ds,~~k=1,2,\cdots,
 \end{equation}all $X_k(\zeta, \mu)$ and $Y_k(\zeta, \mu)$ being analytic functions of $ \zeta$ for finite $\zeta$.
\end{lem}

\noindent
{\bf Proof. }
In view of the Airy equation, by applying Taylor expansions we are capable of writing $\chi(n^{1/3} \zeta +n^{-1/6} \mu)$ as in \eqref{chi-Taylor},
  where the coefficients $X(n; \zeta, \mu)$ and $Y(n; \zeta, \mu)$ are analytic in $\zeta$ and $\mu$, with
 \begin{equation*}
   X(n; \zeta, 0)=1,~~  Y(n; \zeta, 0)=0,
 \end{equation*}and, by taking derivative in \eqref{chi-Taylor} with respect to $\mu$,
we have
\begin{equation}\label{chi-prime-Taylor}
\chi'\left (n^{1/3} \zeta +n^{-1/6} \mu\right )=n^{1/6} \chi\left (n^{1/3} \zeta\right   )\frac{\partial}{\partial \mu}X(n; \zeta, \mu)+     \chi'\left (n^{1/3} \zeta\right )\frac{\partial}{\partial \mu}Y(n; \zeta, \mu).
\end{equation} Hence we also have
 \begin{equation*}
   \frac {\partial}{\partial \mu} X(n; \zeta, 0)=0,~~ \frac {\partial}{\partial \mu} Y(n; \zeta, 0)=1.
 \end{equation*}Furthermore, differentiating \eqref{chi-Taylor}  twice with respect to $\mu$, and compare the resulting equality with \eqref{chi-Taylor}, we obtain the equations
  \begin{equation*}
  \frac {\partial^2}{\partial\mu^2} X(n; \zeta, \mu)=\left (\zeta+\frac {\mu} {n^{1/2}}\right ) X(n; \zeta, \mu), ~~
  X(n; \zeta, 0)=1,~~\frac {\partial}{\partial \mu} X(n; \zeta, 0)=0,
  \end{equation*}and
   \begin{equation*}
  \frac {\partial^2}{\partial\mu^2} Y(n; \zeta, \mu)=\left (\zeta+\frac {\mu} {n^{1/2}}\right ) Y(n; \zeta, \mu), ~~
  Y(n; \zeta, 0)=0,~~\frac {\partial}{\partial \mu} Y(n; \zeta, 0)=1.
  \end{equation*}
 These differential equations, with initial data, have formal (asymptotic) solutions \eqref{X-Y-series},
  where $X_k$ and $Y_k$ are determined, iteratively, by equations
  \begin{equation*}
 \left\{\begin{array}{llll}
         \displaystyle{ \frac {\partial^2 X_0(\zeta, \mu)}{\partial\mu^2}  = \zeta X_0(\zeta, \mu),} & X_0(\zeta, 0)=1,  & \frac {\partial X_0(\zeta, 0) }{\partial \mu}  =0,  &   \\
       \displaystyle{ \frac {\partial^2 X_k(\zeta, \mu)}{\partial\mu^2}  = \zeta X_k(\zeta, \mu)+\mu X_{k-1}(\zeta, \mu)}, & X_k(\zeta, 0)=0 & \frac {\partial X_k(\zeta, 0) }{\partial \mu}  =0, & k=1,2,\cdots,
        \end{array}
 \right .
  \end{equation*}and
  \begin{equation*}
 \left\{\begin{array}{llll}
         \displaystyle{ \frac {\partial^2 Y_0(\zeta, \mu)}{\partial\mu^2}  = \zeta Y_0(\zeta, \mu),} & Y_0(\zeta, 0)=0,  & \frac {\partial Y_0(\zeta, 0) }{\partial \mu}  =1,  &   \\
       \displaystyle{ \frac {\partial^2 Y_k(\zeta, \mu)}{\partial\mu^2}  = \zeta Y_k(\zeta, \mu)+\mu Y_{k-1}(\zeta, \mu)}, & Y_k(\zeta, 0)=0 & \frac {\partial Y_k(\zeta, 0) }{\partial \mu}  =0, & k=1,2,\cdots.
        \end{array}
 \right .
  \end{equation*}
  Solving these equations gives \eqref{X-0-Y-0-sol}, \eqref{X-k-sol} and \eqref{Y-k-sol}. This completes the proof of the Lemma.\qed
  \vskip .5cm

\subsection{Turning point $t=2\sqrt a$ and proof of Theorem \ref{thm-right-turning-point}} \label{subsec:turning-point-right}

 For convenience, we may write for short  $\zeta(t)=\zeta(n, t)=\eta(t)+\Phi(t)/\sqrt n$ and $\mu=u_\pm$. To evaluate $Q_{n\pm 1}(x)$, we need to work out $\chi\left ((n\pm 1)^{1/3}\eta(t_\pm) +(n\pm1)^{-1/6}\Phi(t_\pm)\right )$, that is, $\chi\left(n^{1/3} \zeta +n^{-1/6} u_\pm \right )$; cf.   \eqref{chi-variable-shift}.
  Applying Lemma \ref{lem:chi}, we have
 \begin{equation}\label{chi-Taylor-shift}
 \chi\left(n^{1/3} \zeta +n^{-1/6} u_\pm \right )\sim
 \chi\left(n^{1/3} \zeta \right )\sum^\infty_{k=0} \frac {X_k(\zeta, u_\pm )} {n^{k/2}}+ n^{-1/6}   \chi'\left(n^{1/3} \zeta\right )\sum^\infty_{k=0} \frac {Y_k(\zeta, u_\pm)} {n^{k/2}} .
\end{equation} The situation here   slightly differs from Lemma \ref{lem:chi}  in the dependence on $n$ of  $\zeta$ and $u_\pm$ in lower order terms; see \eqref{u+-leading}-\eqref{u--leading}.
Similarly, in view of \eqref{chi-prime-Taylor}, we obtain
\begin{equation}\label{chi-prime-Taylor-shift}
 \chi'\left(n^{1/3} \zeta +n^{-1/6} u_\pm \right )\sim
n^{1/6} \chi\left(n^{1/3} \zeta \right )\sum^\infty_{k=0} \frac {\frac {\partial}{\partial \mu}X_k(\zeta, u_\pm ) } {n^{k/2}}+     \chi'\left(n^{1/3} \zeta\right )\sum^\infty_{k=0} \frac {\frac {\partial}{\partial\mu} Y_k(\zeta, u_\pm)} {n^{k/2}},
\end{equation} where   $X_k(\zeta,\mu)$ and $Y_k(\zeta,\mu)$ are given in \eqref{X-0-Y-0-sol}-\eqref{Y-k-sol}.

Combining \eqref{chi-Taylor-shift}-\eqref{chi-prime-Taylor-shift} with \eqref{airy-type-sol} yields
\begin{equation}\label{Q-shift}
Q_{n\pm 1}(x)\sim A_\pm(n, \eta)\;  \chi\left(n^{1/3} \zeta \right ) + B_\pm(n, \eta)\; n^{-1/6}   \chi'\left(n^{1/3} \zeta\right ),
\end{equation}where the coefficients demonstrate a complicated dependence on $n$ and $t$ (equivalently, on $\eta$), as
\begin{equation*}
A_\pm(n, \eta)=\sum^\infty_{k=0} \frac {X_k(\zeta, u_\pm )} {n^{k/2}}\sum^\infty_{k=0}  \frac {A_k(\eta_\pm)} {(n\pm 1)^{k/2}} + \left (\frac n {n\pm 1}\right )^{1/6}\sum^\infty_{k=0} \frac {\frac {\partial}{\partial \mu}X_k(\zeta, u_\pm ) } {n^{k/2}}  \sum^\infty_{k=1}  \frac {B_k(\eta_\pm)} {(n\pm 1)^{k/2}}
\end{equation*} and
\begin{equation*}
B_\pm(n, \eta)=\sum^\infty_{k=0} \frac {Y_k(\zeta, u_\pm)} {n^{k/2}}\sum^\infty_{k=0}  \frac {A_k(\eta_\pm)} {(n\pm 1)^{k/2}} +\left (\frac n {n\pm 1}\right )^{1/6}  \sum^\infty_{k=0} \frac {\frac {\partial}{\partial\mu} Y_k(\zeta, u_\pm)} {n^{k/2}}\sum^\infty_{k=1}  \frac {B_k(\eta_\pm)} {(n\pm 1)^{k/2}},
\end{equation*}where $\eta_\pm =\eta(t_\pm)$, $t_\pm$ are the shifted variables; see \eqref{x-invariant}. Finally, substituting \eqref{airy-type-sol} and  \eqref{Q-shift}
 into the difference equation \eqref{charlier-difference-eq-symmetric}, equalizing the coefficients of  $\chi$ and $\chi'$, we have
\begin{equation}\label{chi-coefficients}
 A_+(n, \eta)+A_-(n, \eta)\sim (A_n x+B_n) \sum^\infty_{k=0}  \frac {A_k(\eta)} {n^{k/2}}
\end{equation}and
\begin{equation}\label{chi-prime-coefficients}
 B_+(n, \eta)+B_-(n, \eta)\sim (A_n x+B_n) \sum^\infty_{k=1}  \frac {B_k(\eta)} {n^{k/2}}.
\end{equation}
Now we are in a position to determine the mapping $\eta(t)$. Indeed, picking up  the $O(1)$ terms from \eqref{chi-coefficients}
we have the first order equation
\begin{equation}\label{eta-determine}
e^{-\eta'(t)\sqrt{\eta(t)}} + e^{\eta'(t)\sqrt{\eta(t)}}=\frac t {\sqrt a}.
\end{equation}Here we have made use  of the large-$n$ approximations
\begin{equation*}
  A_nx+B_n\sim \frac t {\sqrt a}-\frac {\sqrt a}{\sqrt n}, ~~~A_0(\eta_\pm)\sim A_0(\eta)\mp \frac {A_0'(\eta) \eta'}{\sqrt n} ,
\end{equation*}
\begin{equation*}
X_0(\zeta, u_\pm)=  \frac {\partial}{\partial \mu} Y_0(\zeta, u_\pm) \sim \frac 1 2\left (e^{-\eta' \sqrt{\eta }} + e^{\eta' \sqrt{\eta }}\right )+ \frac {\sqrt \eta}{ 2\sqrt n} \left ( -\frac {\Phi \eta'}{2\eta} \pm u_{\pm,1}\right ) \left (e^{-\eta' \sqrt{\eta }} - e^{\eta' \sqrt{\eta }}\right )
\end{equation*}and
\begin{equation*}
Y_0(\zeta, u_\pm)\sim \pm \frac 1 { 2\sqrt{\eta }}\left (e^{-\eta' \sqrt{\eta }} - e^{\eta' \sqrt{\eta }}\right )\left ( 1-\frac {\Phi}{2\eta \sqrt n}\right ) +  \frac { 1} {2\sqrt n}\left (u_{\pm, 1}\mp  \frac {\Phi \eta'}{2\eta }\right )\left ( e^{\eta' \sqrt{\eta }}+   e^{-\eta' \sqrt{\eta }} \right ),
\end{equation*} where $u_{\pm, 1}$ are given in \eqref{u+-leading}-\eqref{u--leading}, and
each of the asymptotic approximations has an error $O(1/n)$; cf. \eqref{u+-leading}, \eqref{u--leading}, \eqref{X-0-Y-0-sol} and \eqref{x-invariant}.
While   the leading $O(1)$ terms on both sides of  \eqref{chi-prime-coefficients} vanish.

Near the turning point $t_0$, $\eta(t)$ is uniquely determined from \eqref{eta-determine} by the initial condition $\eta(t_0)=0$, and the assumption that $\eta(t)$ is monotone increasing  for real $t>t_0=2\sqrt a$. Accordingly, we have
\begin{equation}\label{conformal-mapping-derivative}
\eta'(t)\sqrt{\eta(t)}=\log \frac {t+\sqrt{t^2-4a}}{2\sqrt a},~~t\in \mathbb{C}\setminus  (-\infty, 2\sqrt a],
\end{equation}where the logarithm takes principal branch, and  $\sqrt{t^2-4a}$ is analytic in $\mathbb{C}\setminus [-2\sqrt a, 2\sqrt a]$ and is positive for $t>2\sqrt a$. As a result,
\begin{equation}\label{conformal-mapping}
\frac 2 3 \left ( \eta(t) \right )^{3/2} =t\log\frac {t+\sqrt{t^2-4a}}{2\sqrt a}-\sqrt{t^2-4a},~~t\in \mathbb{C}\setminus (-\infty, 2\sqrt a].
\end{equation}It is readily verified that $\frac 2 3 \left ( \eta(t) \right )^{3/2}$ takes purely imaginary values   as $t$ approaches  $(-2\sqrt a, 2\sqrt a)$ from above and below, and that if we choose the branch of $\eta(t)$ in \eqref{conformal-mapping},   to be positive on $t>2\sqrt a$, then $\eta(t)$ is analytic and univalent in $\mathbb{C}\setminus (-\infty, -2\sqrt a\; ]$.

Having had $\eta(t)$, we proceed to determine $\Phi(t)$ and $A_0(\eta)$. To this aim, we  need to work out more details. For example, from \eqref{X-0-Y-0-sol}-\eqref{Y-k-sol}, for free variables $\zeta$ and $\mu$,   we have
\begin{equation*}
X_1(\zeta, \mu)=\frac 1 {8\sqrt\zeta} \left (\mu^2+\frac 1 \zeta\right )\left (e^{\sqrt\zeta \mu}-e^{-\sqrt\zeta \mu}\right )-\frac \mu{8\zeta} \left (e^{\sqrt\zeta \mu}+e^{-\sqrt\zeta \mu}\right )
\end{equation*}and
\begin{equation*}
Y_1(\zeta, \mu)=\frac {\mu^2} {8 \zeta}  \left (e^{\sqrt\zeta \mu}+e^{-\sqrt\zeta \mu}\right )-\frac \mu{8\zeta^{3/2}} \left (e^{\sqrt\zeta \mu}-e^{-\sqrt\zeta \mu}\right ).
\end{equation*}Resuming that $\zeta=\eta(t)+\frac {\Phi(t)}{\sqrt n}$, we have, with errors  $O(1/\sqrt n)$,
\begin{equation*}
X_1(\zeta, u_\pm)\sim X_1(\eta, \mp \eta')=\mp X_1(\eta,   \eta') ,~~Y_1(\zeta, u_\pm)\sim Y_1(\eta, \mp \eta')= Y_1(\eta,  \eta') ,
\end{equation*}  and
\begin{equation*}
\frac{\partial}{\partial\mu} X_0(\zeta, u_\pm)\sim \pm \frac {\sqrt{\eta }} 2\left (e^{-\eta' \sqrt{\eta }} - e^{\eta' \sqrt{\eta }}\right ),~~A_1(\eta_\pm)\sim A_1(\eta),~~B_1(\eta_\pm)\sim B_1(\eta)  .
\end{equation*}
Making use of all these, bring together  the $O(1/\sqrt n)$ terms in \eqref{chi-coefficients} and \eqref{chi-prime-coefficients}, we have
\begin{equation}\label{Phi-diff-eq}
\sqrt\eta  \left( \frac {\Phi\eta'}{2\eta}+\Phi' -\frac \eta 3+\frac {t\eta'} 2\right ) =- \frac  a {\sqrt{t^2-4a}}
\end{equation}and
\begin{equation}\label{A-0-diff}
\frac {\eta'}{\sqrt{\eta}} \frac{\sqrt{t^2-4a}} {\sqrt a}\;  \frac{d A_0(\eta)}{d\eta}+\left [ \left (  \eta'' + \frac {\eta'^2}{2 \eta}\right )\frac  t{2\sqrt a}-\frac {\eta'}{4\eta^{3/2} }\frac {\sqrt{t^2-4a}}{\sqrt a}
\right ] A_0(\eta)=0.
\end{equation}
In view of \eqref{conformal-mapping-derivative}, \eqref{conformal-mapping}, and the fact that $\eta(t_0)=0$, we solve the differential equation \eqref{Phi-diff-eq} to give
$\sqrt {\eta(t)} \Phi(t)=-\frac {t\sqrt{t^2-4a}} 4$,
 which is \eqref{Phi-sol},  where $\sqrt {\eta(t)}$ and    $\sqrt{t^2-4a}$ are positive for $t>t_0=2\sqrt a$.

One can solve $A_0(\eta)$ up to a constant factor. Indeed, we can write \eqref{A-0-diff} as
\begin{equation*}
\frac {dA_0} {A_0}=\left( -\frac t  {2(t^2-4a) \eta'}+\frac 1 {4\eta}\right ) d\eta,
\end{equation*}where $\frac 1 {\eta'} =\frac {d t}{d\eta}$, and $t$ is regarded as a function of $\eta$. Thus $A_0(\eta)$ is determined up to a constant factor independent of both $n$ and $\eta$. From the above equation we readily  pick one solution $A_0(\eta)= \left ( \frac {t^2-4a} {4a\eta}\right )^{-1/4}$, which is \eqref{A-0-sol}.

Now we choose $\chi$ to be $\Ai$ and $\Bi$, and have the following asymptotic formula
\begin{equation}\label{turn-point-right}
C_n^{(a)}(x)\sim (2a)^{\frac n 2} \frac {\Gamma\left((n+1)/ 2\right )} {\Gamma\left( 1/ 2\right )}\left (   K_1(x) Q_n(x)+K_2(x) \tilde Q_n(x)\right ),
\end{equation}
where the asymptotic solutions
\begin{equation*}
Q_n(x)\sim {A_0(\eta)} \Ai \left (n^{1/3}\eta(t) +n^{-1/6}\Phi(t)\right ) ,~~\tilde Q_n(x)\sim {A_0(\eta)} \Bi \left (n^{1/3}\eta(t) +n^{-1/6}\Phi(t)\right ).
\end{equation*}
To determine the coefficients $ K_1(x)$ and $K_2(x)$, again we apply a matching process: we match the approximation \eqref{turn-point-right} with the intermediate asymptotic formula \eqref{approx-intermediate}.
Since $\eta(t)$ is monotone increasing for $t> 2\sqrt a$, and the solution in  \eqref{approx-intermediate} is exponentially small for $t>  2\sqrt a$, we must have $K_2(x)\equiv 0$, that is, the exponentially large part in \eqref{turn-point-right} vanishes. To find $K_1(x)$, we apply  $\Ai(s)\sim \frac {s^{-1/4}}{2\sqrt\pi} e^{-\frac 2 3s^{3/2}}$, $s\to\infty$; cf. \cite[Eq.(9.7.5)]{nist}.  In view of  \eqref{Phi-sol} and \eqref{A-0-sol}, we deduce from \eqref{turn-point-right} that
\begin{equation*}
C_n^{(a)}(x)\sim (2a)^{\frac n 2} \frac {\Gamma\left((n+1)/ 2\right )} {\Gamma\left( 1/ 2\right )}\frac {K_1(x)(4a)^{\frac 1 4}n^{-\frac 1 {12}}(t^2-4a)^{-\frac 1 4}}{2\sqrt\pi}
e^{-\frac 2 3\sqrt n\eta^{\frac 3 2} +\frac {t\sqrt{t^2-4a}} 4},
\end{equation*}where $\frac 2 3 \eta^{ 3/ 2}$ is explicitly given in \eqref{conformal-mapping}.
The approximation will agree with \eqref{approx-intermediate} if
\begin{equation*}
\frac {K_1(x)(4a)^{\frac 1 4}n^{-\frac 1 {12}} }{2\sqrt\pi}\sim  \frac {2^{-\frac  3  4} \pi^{- \frac 1 4}
e^{ \frac a  2}} {\sqrt {w(x)}}.
\end{equation*}
Hence we choose
\begin{equation*}
 K_1(x)=\left ( \frac \pi{2a}\right )^{1/4}e^{\frac a 2} x^{ \frac 1 {12}} w(x)^{-\frac 1 2}.
\end{equation*}Here use has been made of the fact that $x^{\frac 1 {12}}\sim n^{\frac 1 {12}}$ for $n$ large and $|t|\ll \sqrt n$. Substituting it and  $K_2(x)\equiv 0$ into \eqref{turn-point-right}, we obtain the uniform asymptotic approximation \eqref{approximation-turn-point-right} in a neighborhood of the turning point $t=t_0=2\sqrt a$.
\begin{rem}\label{rem:subdominant-sector}
In \eqref{turn-point-right}, $\tilde Q_n(x)$ is dominant in $|\arg\eta|<\frac \pi 3$, and is subdominant in $\frac \pi 3<|\arg\eta|< \pi$, as compared with  $Q_n(x)$. Here
$\arg \eta\sim \arg (t-2\sqrt a\; )$ for small $t-2\sqrt a$. In the matching process we see that the intermediate asymptotics matches the subdominant term involving  $Q_n(x)$ for $|\arg\eta|<\frac \pi 3$, thus $K_2(x)$ must be asymptotically zero, and $K_1(x)$ determined accordingly. Beyond this sector, $K_1(x)$ and $K_2(x)$ are preserved, for $\tilde Q_n(x)$ is subdominant and can not be observed. The ray $t>2\sqrt a$ is a Stokes line. The sector is illustrated in Figure \ref{fig:domains}.
\end{rem}

As a coherence check,  we may specify $t=2\sqrt a\cos\theta\in (-2\sqrt a, 2 \sqrt a)$, $\theta\in (0, \pi)$.
 Recall that $\Ai(-x) \sim \frac 1 {\sqrt\pi x^{1/4}}\cos(\frac 2 3 x^{3/2}-\frac \pi 4)$ as $x\to \infty$ with $|\arg x|<\frac{2\pi} 3$, cf. \cite[Eq.(9.7.9)]{nist}.
For fixed $\theta\in (0, \pi)$, we can rewrite \eqref{approximation-turn-point-right} as
\begin{equation*}
C_n^{(a)}(x)\sim \frac {2C}{\sqrt{w(x)}\; (4a-t^2)^{1/4}}  \cos \left ( \frac  2 3 \sqrt n \left (-\eta \right )^{3/2}+  (-\eta)^{3/2} \frac \Phi \eta-\frac \pi 4   \right ) ,
\end{equation*}where $C$ is the same as in \eqref{approx-between-via-intermediate}, and we have used $x^{\frac 1{12}}\sim n^{\frac 1 {12}}$ in deriving the approximation.
From \eqref{conformal-mapping} and \eqref{Phi-sol}, it is easily seen that
 \begin{equation*}
 \frac  2 3 \left (-\eta \right )^{3/2}=2\sqrt a (\sin\theta-\theta\cos\theta)~~\mbox{and}~~(-\eta)^{3/2} \frac \Phi \eta=\frac a 2 \sin(2\theta).
 \end{equation*} Hence the approximation is exactly  \eqref{approx-between-via-intermediate}, derived from  intermediate asymptotics.

\subsection{Turning point $t=-2\sqrt a$ and proof of Theorem  \ref{thm-left-turning-point}} \label{subsec:turning-point-left}
In this case we employ a slightly different canonical form.
Substituting
 \begin{equation*}
 C_n^{(a)}(x)=(-1)^n (2a)^{\frac n 2} \frac {\Gamma\left((n+1)/ 2\right )} {\Gamma\left( 1/ 2\right )}\mathcal{Q}_n(x)
 \end{equation*}
 into \eqref{charlier-difference-eq} gives the following difference equation
\begin{equation}\label{charlier-difference-eq-symmetric-minus}
\mathcal{Q}_{n+1}(x)+\left (A_n x+B_n\right ) \mathcal{Q}_n(x)+\mathcal{Q}_{n-1}(x)=0,
\end{equation}
where the coefficients $A_n$ and $B_n$ are the same as in
\eqref{difference-eq-balance}.
Once again, we assume that  we have an asymptotic solution to \eqref{charlier-difference-eq-symmetric-minus}  of the form \eqref{airy-type-sol}, and proceed to  determine the functions $\tilde\eta(t)$, $\tilde\Phi(t)$ and $\mathcal{A}_0(\tilde\eta)$ in this case.
Here we use the notations $\tilde\eta(t)$ and $\tilde\Phi(t)$, instead of $\eta(t)$ and $\Phi(t)$  in the previous subsection.

All the derivation leading to \eqref{chi-coefficients}-\eqref{chi-prime-coefficients} holds, one need only to replace  the factor $A_nx+B_n$ with $-(A_nx+B_n)\sim -\frac t {\sqrt a}+\frac {\sqrt a} n$ on the righthand sides. The equation \eqref{eta-determine} now reads
\begin{equation*}
e^{-\tilde\eta'(t)\sqrt{\tilde\eta(t)}} + e^{\tilde\eta'(t)\sqrt{\tilde\eta(t)}}=-\frac t {\sqrt a},
\end{equation*}
where we further require that $\tilde\eta(-2\sqrt a )=0$, and $\tilde\eta(t)$ is monotone decreasing for $t<-2\sqrt a$. Hence we have
\begin{equation}\label{conformal-mapping-minus}
\frac 2 3 \tilde\eta^{3/2} =t\log \frac{-t+\sqrt{t^2-4a}}{2\sqrt a}+\sqrt{t^2-4a},
\end{equation}such that $\tilde\eta(t)>0$ for $t<-2\sqrt a$,
 ${\tilde\eta}(t)\sim   a^{-1/6} (-2\sqrt a-t )$
for $t\sim -2\sqrt a$, and $\tilde\eta' \sqrt {\tilde\eta}=\log \frac {-t+\sqrt{t^2-4a}}{2\sqrt a}$.

Now compare the $O(1/\sqrt n)$
terms in the modified version  \eqref{chi-coefficients}-\eqref{chi-prime-coefficients}, instead of \eqref{Phi-diff-eq} and \eqref{A-0-diff}, we have
\begin{equation*}
\sqrt{\tilde\eta } \left( \frac {{\tilde\Phi}\tilde\eta'}{2\tilde\eta}+{\tilde\Phi}' -\frac {\tilde\eta} 3+\frac {t\tilde\eta'} 2\right ) =  \frac  a {\sqrt{t^2-4a}}
\end{equation*}and
\begin{equation*}
\frac {\tilde\eta'}{\sqrt{\tilde\eta}} \frac{\sqrt{t^2-4a}} {\sqrt a}\;  \frac{d \mathcal{A}_0(\tilde\eta)}{d\tilde\eta}+\left [ \left (  \tilde\eta'' + \frac {\tilde\eta'^2}{2 \tilde\eta}\right )\frac  {-t}{2\sqrt a}-\frac {\tilde\eta'}{4\tilde\eta^{3/2} }\frac {\sqrt{t^2-4a}}{\sqrt a}
\right ] \mathcal{A}_0(\tilde\eta)=0.
\end{equation*}
The first equation gives $(\sqrt {\tilde\eta} {\tilde\Phi})'=\frac {\sqrt{t^2-4a}} 2+\frac a {\sqrt{t^2-4a}}$, and in turn gives
$
 \sqrt {\tilde\eta} {\tilde\Phi} =\frac  1 4 t \sqrt{t^2-4a}$.

 To solve $\mathcal{A}_0(\tilde\eta)$ from the second equation, in view of \eqref{conformal-mapping-minus}, we can write
\begin{equation*}
\frac {d\mathcal{A}_0} {\mathcal{A}_0}=\left( -\frac t  {2(t^2-4a) \tilde\eta'}+\frac 1 {4\tilde\eta}\right ) d\tilde\eta,
\end{equation*}the same as in the previous subsection.
Hence  $\mathcal{A}_0(\tilde\eta)$ is determined up to a constant factor independent of both $n$ and $\tilde\eta$. From the above equation we readily  pick one solution $
\mathcal{A}_0(\tilde\eta)= \left ( \frac {t^2-4a} {4a\tilde\eta}\right )^{-1/4}$ in the sense of \eqref{A-0-sol-left}.
 It is worth pointing out that $\mathcal{A}_0(\tilde\eta)$ is, like $\tilde\eta(t)$,   analytic at $t=-2\sqrt a$, such that
$\mathcal{A}_0(\tilde\eta)\sim a^{1/8} \left (\frac {-(t+2\sqrt a)}{\tilde\eta}\right )^{-1/4}$ as $t\to-2\sqrt a$, and  is real positive for $t<-2\sqrt a$, where
$\arg \{-(t+2\sqrt a)\} \in (-\pi, \pi)$.

Denote by $\mathcal{Q}_n(x)$ and $\tilde{\mathcal{Q}}_n(x)$ the asymptotic solutions to \eqref{charlier-difference-eq-symmetric-minus}, such that
\begin{equation*}
\mathcal{Q}_n(x)\sim {\mathcal{A}_0(\tilde\eta)} \Ai \left (\omega n^{\frac 1 3}\left (\tilde\eta  +\frac {{\tilde\Phi}}{\sqrt n}\right )\right ) ,~~\tilde {\mathcal{Q}}_n(x)\sim {\mathcal{A}_0(\tilde\eta)} \Ai \left (\omega^2  n^{\frac 1 3}\left (\tilde\eta  +\frac {{\tilde\Phi} }{\sqrt n}\right )\right ),
\end{equation*}where $\omega =e^{ 2\pi i/3 }$,  functions $\tilde\eta(t)$, ${\tilde\Phi}(t)$ and $\mathcal{A}_0(\tilde\eta)$ are given in the present subsection.  Accordingly we can write
\begin{equation}\label{turn-point-left}
C_n^{(a)}(x)\sim (-1)^n (2a)^{\frac n 2} \frac {\Gamma\left((n+1)/ 2\right )} {\Gamma\left( 1/ 2\right )}\left (   \mathcal K_1(x) \mathcal Q_n(x)+\mathcal K_2(x) \tilde {\mathcal Q}_n(x)\right ),
\end{equation}with $\mathcal K_1(x)$ and $\mathcal K_2(x)$ to be determined by matching \eqref{turn-point-left} with the intermediate asymptotic formula \eqref{approx-intermediate}.
As in Remark \ref{rem:subdominant-sector}, we pay attention to the dominant and subdominant solutions. First we consider the case when $\arg \tilde\eta \sim -\frac {2\pi} 3$, or, approximately,
$\arg(t+2\sqrt a)= \pi +\arg (-(t+2\sqrt a))\sim  \frac {\pi} 3$. Actually, the curve $\arg \tilde\eta(t) = -\frac {2\pi} 3$ is a Stokes line; cf. Figure \ref{fig:domains}.    In this case,
$\tilde {\mathcal Q}_n(x)$ is dominant as compared with $\mathcal Q_n(x)$.
The  dominant  $\tilde {\mathcal Q}_n(x)$ does not match \eqref{approx-intermediate}, hence we set $\mathcal K_2(x)=0$. To determine $\mathcal K_1(x)$, we expand \eqref{turn-point-left} to give
\begin{equation*}
\frac {C_n^{(a)}(x)} {(2a)^{\frac n 2} \frac {\Gamma\left((n+1)/ 2\right )} {\Gamma\left( 1/ 2\right )}}   \sim     \mathcal K_1(x)\frac {e^{\frac {\pi i} 3} (4a)^{\frac 1 4}n^{-\frac 1 {12}}}{2\sqrt\pi} e^{x\pi i+\sqrt n\left[ t\log\frac {t-\sqrt{t^2-4a}}{2\sqrt a}+\sqrt{t^2-4a}\right ]-\frac 1 4\log(t^2-4a)+\frac 1 4t\sqrt{t^2-4a}}.
\end{equation*}Comparing  it with \eqref{approx-intermediate}, and recalling that
$x^{\frac 1 {12}}\sim n^{\frac 1 {12}}$,
we obtain
\begin{equation*}
\mathcal K_1(x)=\left (\frac \pi {2a}\right )^{1/4} e^{a/2} x^{\frac 1 {12}} (w(x))^{-1/2} e^{-(x\pi+\frac \pi 3)i}.
\end{equation*}Indeed, the above matching process holds for  $\arg \tilde\eta \in (-\pi, -\frac \pi 3)$, in which the intermediate asymptotics matches the subdominant solution $\mathcal Q_n(x)$.
For $\arg \tilde\eta \in (-\frac \pi 3, 0)$, the result is still valid since $\mathcal Q_n(x)$ is the dominant solution. Therefore, for $t+2\sqrt a=O(1)$, $\Im t>0$, we have \eqref{approximation-turn-point-left} in the upper half $t$-neighborhood, where the special function employed  is an Airy function $e^{-{\pi i} /3}\Ai \left ( s e^{2\pi i/3}\right )=\frac 1 2\left ( \Ai(s)- i \Bi(s)\right )$; see \cite[Eq.(9.2.11)]{nist}.  The formula \eqref{approximation-turn-point-left-real} follows from taking twice real part of \eqref{approximation-turn-point-left}.
The asymptotic approximation in the lower half neighborhood is obtained by symmetry, involving the Airy function $e^{  {\pi i} /3}\Ai \left ( s e^{- 2\pi i/3}\right )=\frac 1 2\left ( \Ai(s)+ i \Bi(s)\right )$.
 This proves Theorem \ref{thm-left-turning-point}.

As an application of Theorem \ref{thm-left-turning-point}, we check the oscillating  in an real interval  around $t=-2\sqrt a$. When $t\sim -2\sqrt a$, the zeros of the Charlier polynomials will be represented in terms of the zeros of Airy function. Here, we would rather do a coherence check to show that how the density of zeros changes from \eqref{approx-in-Charlier-formula} to \eqref{approx-between-via-intermediate}. To this aim, we require $t$ to keep away from $-2\sqrt a$.

For $t<-2\sqrt a$, $s=\tilde\eta(t)+{\tilde\Phi}(t)/\sqrt n$ is positive, and $\Bi(n^{1/3} s)$ is dominant over $\Ai(n^{1/3} s)$. Hence from \eqref{approximation-turn-point-left-real} we deduce
\begin{equation*}
C_n^{(a)}(x)\sim C_{\mathcal{K}, n} \frac { x^{\frac 1 {12}}\mathcal{A}_0(\tilde\eta)}{\sqrt{w(x)} } \Bi \left (n^{\frac 1 3}\left (\tilde\eta  +\frac {{\tilde\Phi}}{\sqrt n}\right )\right )
\cos \left (x\pi +\frac \pi 2\right ).
\end{equation*}Here use has been made of the fact that $\Bi(s)\sim \frac 1 {\sqrt \pi s^{1/4}} e^{\frac 2 3 s^{3/2}}$ for large positive $s$; see \cite[Eq.(9.7.7)]{nist}.
It is readily seen that the above formula agrees with  \eqref{approx-in-Charlier-formula} as $1\ll -(t+2\sqrt a)\ll \sqrt n$.

On the other side of $t=-2\sqrt a$, namely, $t>-2\sqrt a$,   $s=\tilde\eta(t)+{\tilde\Phi}(t)/\sqrt n$ is negative. Parameterizing $t=2\sqrt a \cos\theta$ with $\theta<\pi$, and using the asymptotic formulas \cite[Eq.(9.7.9), Eq.(9.7.11)]{nist}   for $\Ai$ and $\Bi$, from \eqref{approximation-turn-point-left} one deduce
\begin{equation*}
C_n^{(a)}(x)\sim (-1)^n C_{\mathcal{K}, n} \frac {\mathcal{A}_0(\tilde\eta)}{\sqrt{\pi w(x)}(-\tilde\eta)^{1/4}}\;  \cos\left (2\sqrt {a   n } \left [\sin\theta -\theta\cos\theta \right ]+a\sin\theta\cos\theta -\frac \pi 4    \right )
\end{equation*}for $t>-2\sqrt a$, as $n\to\infty$; see the discussion right after \eqref{approx-in-Charlier-formula}.
Here we have used
\begin{equation*}
\frac 2 3(-\tilde\eta)^{3/2}=2\sqrt a \left [ (\pi-\theta)\cos\theta+\sin\theta\right ],~~ (-\tilde\eta)^{3/2}\frac  {\tilde\Phi} {\tilde\eta} =a\cos\theta\sin\theta~~\mbox{and}~~\mathcal{A}_0(\tilde\eta)=\left (\frac {-\tilde\eta}{\sin^2\theta }\right )^{1/4}.
\end{equation*}
The above approximation for $C_n^{(a)}(x)$ is exactly \eqref{approx-between-via-intermediate}.

\section{Comparison with earlier literature}\label{sec:comparison}
We compare  local behavior with   the asymptotic formulas  obtained in Bo and Wong \cite{bo-wong1994}.

The first case is the asymptotic approximation for $C^{(a)}_n(ny)$ for $y\in [\epsilon, r]\subset (0, 1)$. The leading term formula given in \cite[p.310]{bo-wong1994} reads
\begin{equation}\label{Bo-Wong-at-0}
\frac 1 {n!} C^{(a)}_n(ny)\sim \left (\frac y {1-y}\right )^{1/2}
e^{a/(1-y)} y^{ny} (1-y)^{n(1-y)}\sin [n(1-y)\pi]\sqrt {\frac 2 {n\pi}},~~~n\to\infty;
\end{equation}see also Goh \cite[Eq.(84)]{goh1998}. This is the situation when Theorem \ref{thm-approx-not-at-1} applies. Indeed, we can write \eqref{approx-not-at-1} as
\begin{equation*}
C^{(a)}_n(ny)  \sim \frac 1 \pi {(-1)^{n+1}} e^{ \frac {a }{1-y }}  {\Gamma\left (n-ny\right )}{\Gamma(1+ny)}\sin (ny\pi),~~~n\to\infty.
\end{equation*}Noting that $n-ny\geq (1-r)n\gg 1$ and $1+ny>\epsilon n\gg 1$,
and applying Stirling's formula, we obtain \eqref{Bo-Wong-at-0}.
It is worth pointing out that if $\epsilon=0$,
 the approximations \eqref{approx-not-at-1} and \eqref{Bo-Wong-at-0} will no longer agree with each other; cf. the case when $ny=O(1)$.

The next case is when $y$ is near $1+2\sqrt a /\sqrt n$, in a way described in \cite[Eq.(5.3)]{bo-wong1994}, namely
\begin{equation*}
y=1+\frac {2\sqrt a}{\sqrt n} +\frac {s}{n^{5/6}}+\frac a n.
\end{equation*}
Local behavior in this case is
\begin{equation}\label{Bo-Wong-at-1+}
\frac 1 {n!} C^{(a)}_n(ny)\sim e^{\frac{3a} 2}\left (\frac n{a}\right )^{\sqrt{an} +(s/2) n^{1/6} +a/2}  (an)^{-\frac 1 6}  \Ai\left(  s a^{-1/6}\right )~~~\mbox{as}~~n\to\infty;
\end{equation}cf. Bo and Wong \cite[Eq.(5.11)]{bo-wong1994}  and  Goh \cite[Eq.(30)]{goh1998}.
We see that Theorem \ref{thm-right-turning-point} applies to this special case, with $y=1+t/\sqrt n$, and
\begin{equation*}
t-2\sqrt a=\frac s {n^{1/3}}+\frac a {\sqrt n},
\end{equation*}small for finite $s$ and large $n$. From \eqref{conformal-mapping-introduction} and \eqref{Phi-sol}, straightforward calculation gives
\begin{equation*}
\eta(t)=a^{-\frac 1 6}(t-2\sqrt a)+O  \left((t-2\sqrt a)^2\right )
=\frac{a^{-1/6}s}{n^{1/3}}+\frac{a^{5/6} }{\sqrt n }+O \left(\frac 1{n^{2/3}} \right ),
\end{equation*}
and
\begin{equation*}
\Phi(t)=-a^{\frac 5 6}+O \left(\frac 1{n^{1/3}} \right ).
\end{equation*}Hence we have
\begin{equation*}
n^{\frac 1 3}\eta(t) +n^{-\frac 1 6}\Phi(t)=a^{-\frac 1 6}s+O\left(\frac 1 {n^{1/3}}\right ) .
\end{equation*}
Substituting it into \eqref{approximation-turn-point-right}, we have
  \begin{equation*}
  C^{(a)}_n(ny)\sim e^{\frac{3a} 2}\sqrt{2\pi}\; e^{(n+\frac 1 2)\log n-n}\left (\frac n{a}\right )^{\sqrt{an} +(s/2) n^{1/6} +a/2}  (an)^{-\frac 1 6}  \Ai\left(  s a^{-1/6}\right )~~~\mbox{as}~~n\to\infty,
  \end{equation*}
which agrees with
\eqref{Bo-Wong-at-1+}.

Our last comparison will be made in yet another special case
\begin{equation*}
y=1-\frac {2\sqrt a}{\sqrt n} +\frac {s}{n^{5/6}}+\frac a n;
\end{equation*}cf. \cite[p.311]{bo-wong1994}. Bo and Wong checked  this case, so as to compare their asymptotic formula \cite[Eq.(5.12)]{bo-wong1994}, with a result of Goh \cite[Eq.(51)]{goh1998}  in terms of Scorer's function.
The formula of Bo and Wong states that
\begin{equation}\label{Bo-Wong-at-1-}
\frac 1 {n!} C^{(a)}_n(ny)\sim e^{\frac{3a} 2}\left (\frac n{a}\right )^{(x-n)/2} (an)^{-\frac 1 6}(-1)^n 2\Re\left [ e^{(x\pi+\pi/3)i}\Ai\left(e^{\pi i/3} s a^{-1/6}\right )\right ]
\end{equation}as $n\to\infty$, where $x=ny=n-2\sqrt{an}+s n^{1/6}+a$.

With our local transformation $y=1+t/\sqrt n$, this case corresponds to
\begin{equation*}
-2\sqrt a -t=\frac {-s} {n^{1/3}}-\frac a {\sqrt n}.
\end{equation*}We see that $t$ is near the turning point $-2\sqrt a$, and
  Theorem \ref{thm-left-turning-point} applies. From \eqref{conformal-mapping-minus-introduction} and \eqref{Phi-sol-left} we find
  \begin{equation*}
  \tilde\eta(t)  =  a^{-1/6} (-2\sqrt a-t )+O  \left((-2\sqrt a-t )^2\right )
=\frac{a^{-1/6}(-s)}{n^{1/3}}-\frac{a^{5/6} }{\sqrt n }+O \left(\frac 1{n^{2/3}} \right ),
\end{equation*}
and
\begin{equation*}
\tilde\Phi(t)=a^{\frac 5 6}+O \left(\frac 1{n^{1/3}} \right ).
\end{equation*}Hence we have the real variable
\begin{equation*}
\Theta:=n^{\frac 1 3}\left(\tilde \eta(t) +\frac {\tilde\Phi(t)}{\sqrt n}\right )=a^{-\frac 1 6}(-s)+O\left(\frac 1 {n^{1/3}}\right ) .
\end{equation*}
Also, $\mathcal{A}_0({\tilde\eta})$ in
   \eqref{A-0-sol-left}, analytic at $t=-2\sqrt a$, has the leading behavior
   \begin{equation*}
\mathcal{A}_0({\tilde\eta})= a^{\frac 1 8} \left (\frac {-2\sqrt a-t }{\tilde\eta}\right )^{-1/4}+O\left ( -2\sqrt a-t\right )=a^{\frac 1 {12}} +O \left(\frac 1{n^{1/3}} \right ).
\end{equation*}
Substituting these into an equivalent form of \eqref{approximation-turn-point-left-real}, that is,
\begin{equation*}
C_n^{(a)}(x)\sim \frac {C_{\mathcal{K}, n} x^{\frac 1 {12}}\mathcal{A}_0({\tilde\eta})}{\sqrt{w(x)} }
  \left [ e^{(x\pi+\pi/3)i}\Ai\left(e^{-2\pi i/3}\Theta\right )+ e^{-(x\pi+\pi/3)i}\Ai\left(e^{2\pi i/3}\Theta\right )\right ] ;
\end{equation*}see \cite[Eq.(9.2.11)]{nist}, and applying Stirling's formula, we obtain
 \begin{eqnarray*}
  C_n^{(a)}(x) & \sim&  e^{\frac{3a} 2}(-1)^n \sqrt{2\pi} e^{\left(n+\frac 1 2\right)\log n-n} (an)^{-\frac 1 6} \left (\frac n a \right )^{(x-n)/2} \\
     & &  \times \left [ e^{(x\pi+\pi/3)i}\Ai\left(e^{\pi i/3}s a^{-1/6}\right )+ e^{-(x\pi+\pi/3)i}\Ai\left(e^{-\pi i/3}s a^{-1/6}\right )\right ] ,
 \end{eqnarray*}
which agrees with
\eqref{Bo-Wong-at-1-}.

\section*{Acknowledgements}
The authors are very grateful to Dr. Dan Dai, Dr. Yu-Tian Li and Prof. Roderick Wong for their support, and to the anonymous reviewers for their helpful and constructive comments.
 X.-M. Huang was partially supported by grants from the Research Grants Council of the Hong Kong Special Administrative Region, China (Project No. CityU 11300115, CityU 11303016).  
The work of Y. Lin was supported in part by the National Natural Science Foundation of China Under Grant Number 11501215, GuangDong Natural Science Foundation under grant number 2014A030310092 and the Fundamental Research Funds for the Central Universities under grant number 2015ZM090.
Y.-Q. Zhao  was supported in part by the National
Natural Science Foundation of China under grant numbers 11571375 and 11971489.


\begin{thebibliography}{99}


\bibitem{askey1985}R. Askey,  Review of the book  {\it An introduction to orthogonal polynomials}, by T.S. Chihara, Gordon \& Breach, New York/London/Paris, 1978, 
{\it J. Approx. Theory},   {\bf 43} (1985),
 394-395.




\bibitem{baik-et-al2007} J. Baik, T. Kriecherbauer, K.T.-R. McLaughlin and P.D. Miller,
Discrete orthogonal polynomials, asymptotics and applications,
{\it{Ann. Math. Studies}}, {\bf{164}}, Princeton University Press.
Princeton and Oxford, 2007.



\bibitem{bender-orszag-1999-book} C.M. Bender and S.A. Orszag, {\it Advanced mathematical methods for scientists and engineers},
Reprint of the 1978 original, Springer-Verlag, New York, 1999.


\bibitem{bo-wong1994}R. Bo and R. Wong,
Uniform asymptotic expansion of Charlier polynomials,
{\it Methods Appl. Anal.}, {\bf  1} (1994),   294-313.


\bibitem{cao-li2014}
L.-H. Cao and Y.-T. Li, Linear difference equations with a transition point at  the origin,
{\it Anal. Appl.}, {\bf 12} (2014), 75-106.





\bibitem{chihara-book} T.S.  Chihara, {\it An introduction to orthogonal polynomials}, Gordon and Breach, New York, 1978.

\bibitem{costin-costin1996}
O. Costin and R. Costin,  Rigorous WKB for finite-order linear recurrence relations with smooth coefficients, {\it SIAM J. Math. Anal.}, {\bf 27} (1996),   110-134.


\bibitem{dai-ismail-wang}D. Dai,
M.E.H. Ismail and
X.-S. Wang, Plancherel-Rotach asymptotic expansion for some polynomials from indeterminate moment problems,
{\it Constr.
Approx.}, {\bf{40}} (2014), 61-104.

\bibitem{deift1999} P. Deift, {\it{Orthogonal polynomials and random matrices:
a Riemann-Hilbert approach}}, Courant Lecture Notes 3, New York
University, 1999.


\bibitem{dunster2001}
T.M. Dunster,  Uniform asymptotic expansions for Charlier polynomials, {\it J. Approx. Theory},   {\bf 112}  (2001),   93-133.

\bibitem{geronimo2009}
J.S. Geronimo, WKB and turning point theory for second order difference equations:
external fields and strong asymptotics for orthogonal polynomials, {\it
arXiv:0905.1684.}



\bibitem{goh1998}W.M.Y. Goh,   Plancherel-Rotach asymptotics for the Charlier
polynomials,
{\it Constr. Approx.}, {\bf 14} (1998),   151-168.

\bibitem{huang2018}X.-M. Huang, {\it The difference equation method to asymptotic expansion of orthogonal polynomials and confluent Heun equation},
Ph.D  Thesis, Sun Yat-sen University, 2018.

\bibitem{huang-cao-wang2017}X.-M. Huang, L.-H. Cao and X.-S. Wang,  Asymptotic expansion of orthogonal polynomials via difference equations, {\it  J. Approx. Theory}, {\bf{239}} (2019),  29-50.

\bibitem{kuijlaars-vanassche1999} A.B.J. Kuijlaars and W. Van Assche,    The asymptotic zero distribution of orthogonal polynomials with varying recurrence coefficients,  {\it  J. Approx. Theory}, {\bf{99}} (1999),  167-197.

\bibitem{maejima-vanassche1985} M. Maejima and  W. Van Assche, Probabilistic proofs of asymptotic formulas for some classical
polynomials, {\it Math. Proc. Cambridge Philos. Soc.}, {\bf 97} (1985), 499-510.

\bibitem{nist}F. Olver, D. Lozier, R. Boisvert and C. Clark,
{\it NIST handbook of mathematical functions},  Cambridge University Press, Cambridge, 2010.



\bibitem{ou-wong2010} C.-H. Ou and R. Wong,  The Riemann-Hilbert approach to global asymptotics of discrete orthogonal polynomials with infinite nodes, {\it Anal. Appl.}, {\bf 8} (2010), 247-286.




\bibitem{szego-book} G. Szeg\H{o}, {\it Orthogonal polynomials}, Fourth edition, American Mathematical Society, Providence, Rhode Island, 1975.


\bibitem{vanassche-geronimo1989}
W. Van Assche and J.S. Geronimo,
Asymptotics for orthogonal polynomials with regularly varying recurrence coefficients,
{\it
Rocky Mountain J. Math.},  {\bf 19}  (1989), 39-49.



\bibitem{wang-wong2012}X.-S. Wang and  R. Wong, Asymptotics of orthogonal polynomials via recurrence relations,
{\it
Anal. Appl.}, {\bf{10}} (2012), 215-235.

\bibitem{wang-wong2002}
Z. Wang  and R. Wong,
Uniform asymptotic expansion of $J_\nu(\nu a)$ via a difference equation, {\it Numer. Math.}, {\bf 91} (2002), 147-193.

\bibitem{wang-wong2003}
Z. Wang and R. Wong, Asymptotic expansions for second-order linear difference equations with a turning  point,
{\it Numer. Math.}, {\bf 94} (2003), 147-194.

\bibitem{wang-wong2005}
Z. Wang and  R. Wong,  Linear difference equations with transition points,
{\it Math. Comp.}, {\bf 74} (2005), 629-653.

\bibitem{wong2014}R. Wong, Asymptotics of linear recurrences,
{\it Anal. Appl.}, {\bf 12} (2014), 463-484.


\bibitem{wong-li-1992}R. Wong and H. Li, Asymptotic expansions for second-order linear difference equations,
{\it  J. Comput. Appl. Math.}, {\bf 41} (1992),  65-94.









\end{thebibliography}
\end{document}